\documentclass{qam-l}

\usepackage{amsmath,amssymb}
\usepackage{graphicx,color}
\usepackage{float}

\usepackage{amsfonts,graphicx}

\renewcommand{\Re}{\mathbb{R}}

\newcommand{\half}{\frac{1}{2}}
\newcommand{\weak}{\rightharpoonup}
\newcommand{\weakst}{\rightharpoonup^*}
\newcommand{\embed}{\hookrightarrow}

\newcommand{\vph}{\vphantom{A^{A}_{A}}}
\newcommand{\sst}{\,\mid\,}

\newcommand{\dv}{\mathrm{div}}

\newcommand{\eqnref}[1]{(\ref{eqn#1})}

\newcommand{\uhat}{\hat{u}}

\newcommand{\gammahat}{\hat{\gamma}}

\newcommand{\etahat}{\hat{\eta}}

\newcommand{\nuhat}{\hat{\nu}}

\newcommand{\phihat}{\hat{\phi}}

\newcommand{\bbC}{\mathbb{C}}

\newcommand{\bbG}{\mathbb{G}}

\newcommand{\bbN}{\mathbb{N}}

\newcommand{\bbU}{\mathbb{U}}

\newcommand{\ubar}{{\bar{u}}}

\newcommand{\Sbar}{{\bar{S}}}
\newcommand{\Tbar}{{\bar{T}}}

\newcommand{\gammabar}{\bar{\gamma}}

\newcommand{\ft}{\tilde{f}}

\newcommand{\ut}{\tilde{u}}

\newcommand{\calA}{{\cal A}}
\newcommand{\calB}{{\cal B}}

\newcommand{\calL}{{\cal L}}

\newcommand{\calU}{{\cal U}}

\newcommand{\norm}[1]{\| {#1} \|}
\newcommand{\Hone}{{H^1(\Omega)}}

\newcommand{\Ltwo}{{L^2(\Omega)}}
\newcommand{\ltwo}[1]{\norm{#1}_\Ltwo}

\newcommand{\Lp}{{L^p(\Omega)}}

\newcommand{\Wonefour}{{W^{1,4}(\Omega)}}

\newcommand{\LoneLtwo}{{L^1[0,T;\Ltwo]}}
\newcommand{\loneltwo}[1]{\norm{#1}_\LoneLtwo}

\newcommand{\LinfLtwo}{{L^\infty[0,T;\Ltwo]}}
\newcommand{\linfltwo}[1]{\norm{#1}_\LinfLtwo}

\newcommand{\LtwoHone}{{L^2[0,T;\Hone]}}

\newcommand{\CHone}{{C[0,T;\Hone]}}

\newcommand{\Curl}{\mathrm{Curl}}

\newcommand{\cal}{\mathcal}

\newtheorem{theorem}{Theorem}[section]
\newtheorem{lemma}[theorem]{Lemma}
\newtheorem{corollary}[theorem]{Corollary}

\theoremstyle{definition}

\newtheorem{notation}[theorem]{Notation}

\begin{document}

\bibliographystyle{amsplain}

\title{Analysis of a Dislocation Model for Earthquakes\footnote{Submitted
Quarterly of Applied Mathematics, April 2016.}}

\author{Jing Liu}
\address{Department of Mathematical Sciences, Carnegie
    Mellon University, Pittsburgh, PA 15213.}
\email{jingliu1@@andrew.cmu.edu}

\author{Xin Yang Lu}
\address{Department of Mathematics and Statistics, 
  McGill University, Montreal, Canada}
\email{xinyang.lu@@mcgill.ca}

\author{Noel J.  Walkington}
\address{Department of Mathematical Sciences, Carnegie
    Mellon University, Pittsburgh, PA 15213.}
\email{noelw@@andrew.cmu.edu}
\thanks{Supported in part by National Science
    Foundation grants DMS--1418991 and DMREF--1434734. This work was
    also supported by the NSF through the Center for Nonlinear
    Analysis.}


\date{\today}

\begin{abstract}
  Approximation of problems in linear elasticity having small shear
  modulus in a thin region is considered. Problems of this type arise
  when modeling ground motion due to earthquakes where rupture occurs
  in a thin fault. It is shown that, under appropriate scaling,
  solutions of these problems can be approximated by solutions of a
  limit problem where the fault region is represented by a surface.
  In a numerical context this eliminates the need to resolve the large
  deformations in the fault; a numerical example is presented to
  illustrate efficacy of this strategy.
\end{abstract}

\maketitle

%\begin{keywords}
%Gamma convergence, earthquake simulation, numerical approximation.
%\end{keywords}

\section{Introduction}\label{sec:intro}
\begin{figure}
\hspace{-0.8cm}
\includegraphics[height=2.25in,clip=true, bb=144 80 630 480]{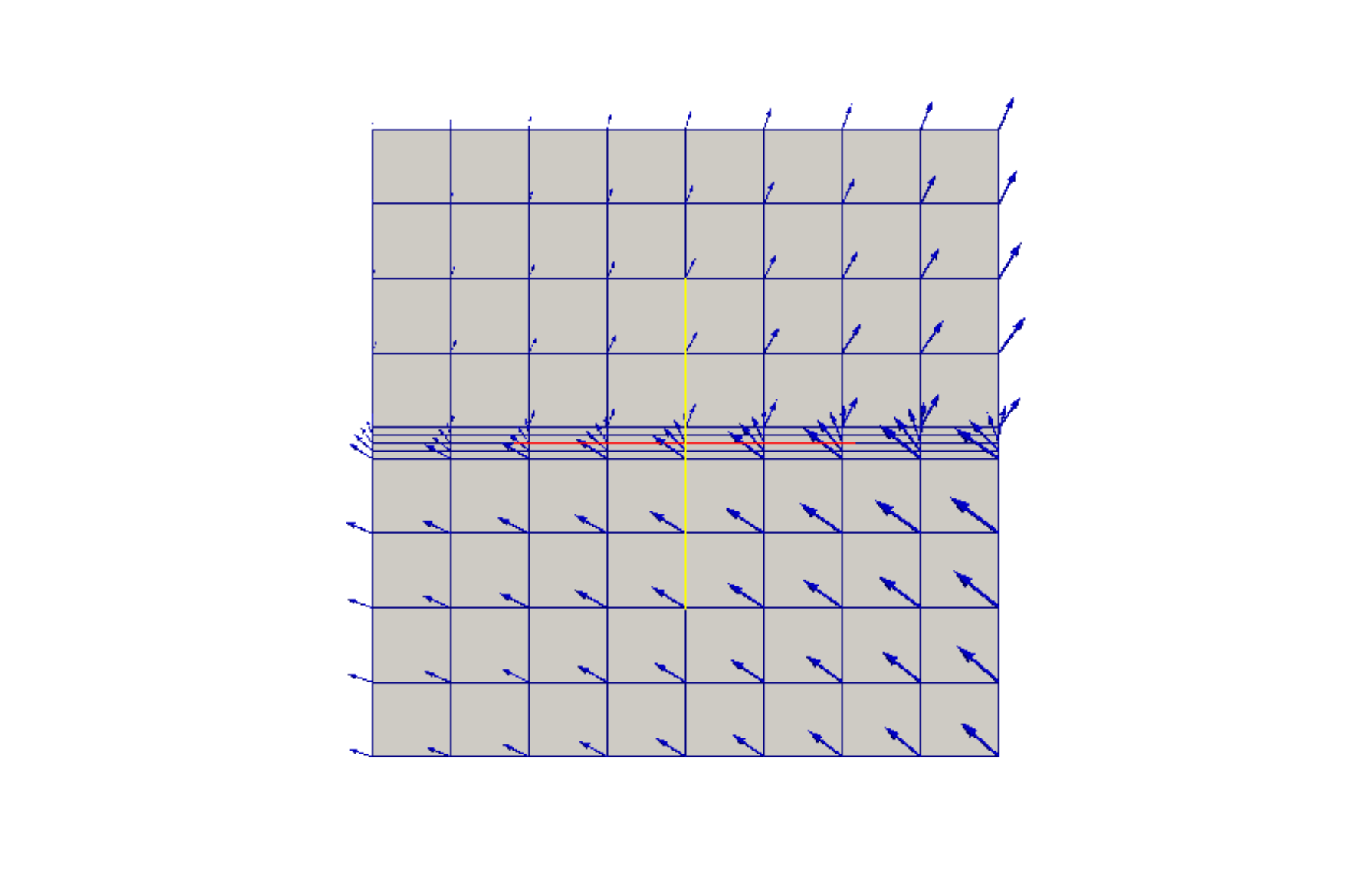} 
\hfill
\includegraphics[height=2.25in,clip=true, bb=155 520 440 785]{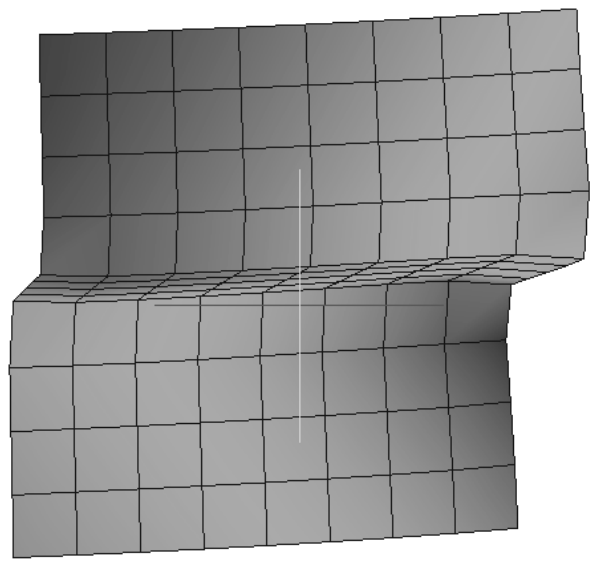}
\caption{Solution with deformation resolved.} 
\label{fig:elastic3m844}
\end{figure}

\begin{figure}
\includegraphics[height=2.25in,clip=true, bb=144 80 600 480]{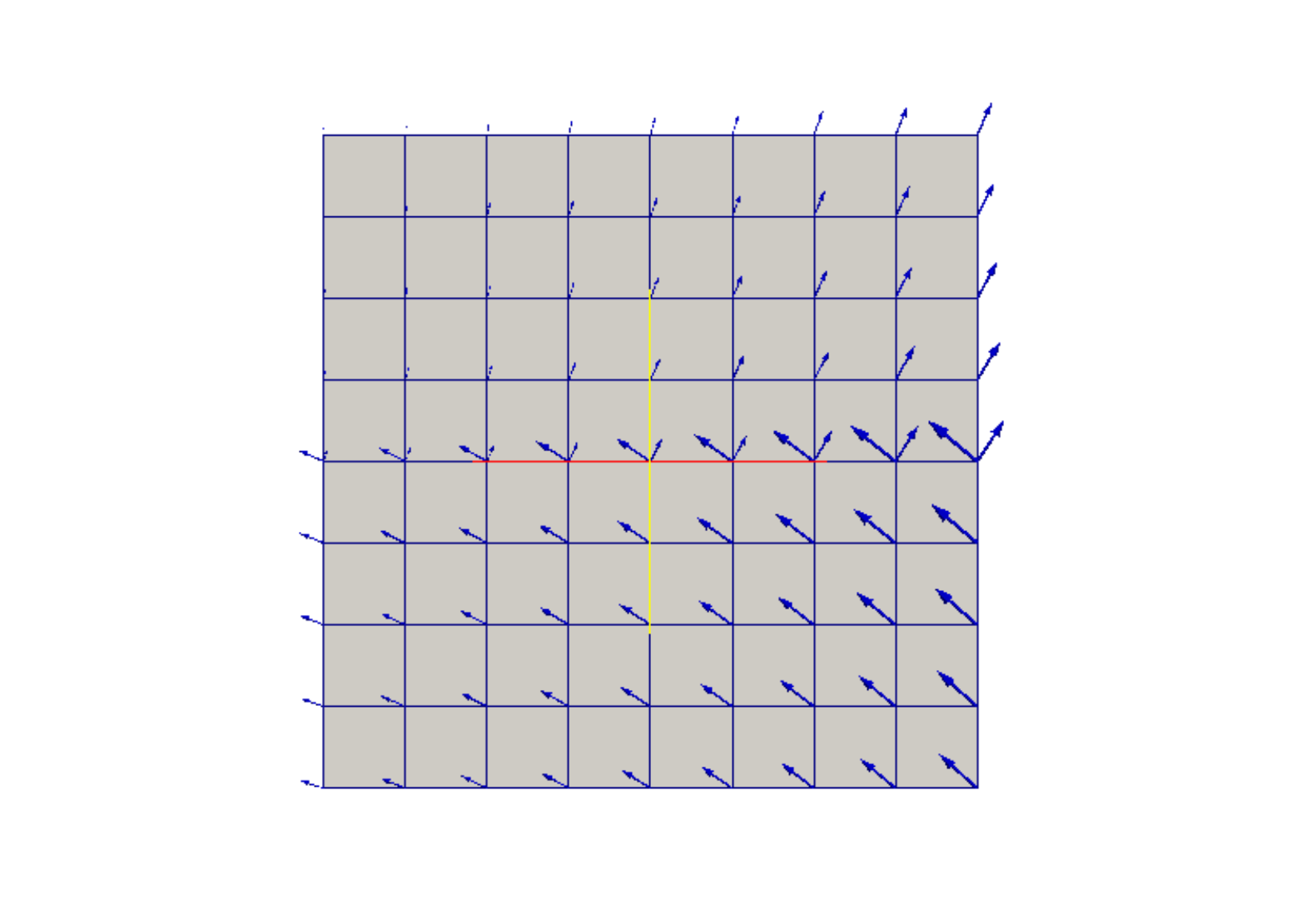}
\hfill
\includegraphics[height=2.25in,clip=true, bb=150 485 450 775]{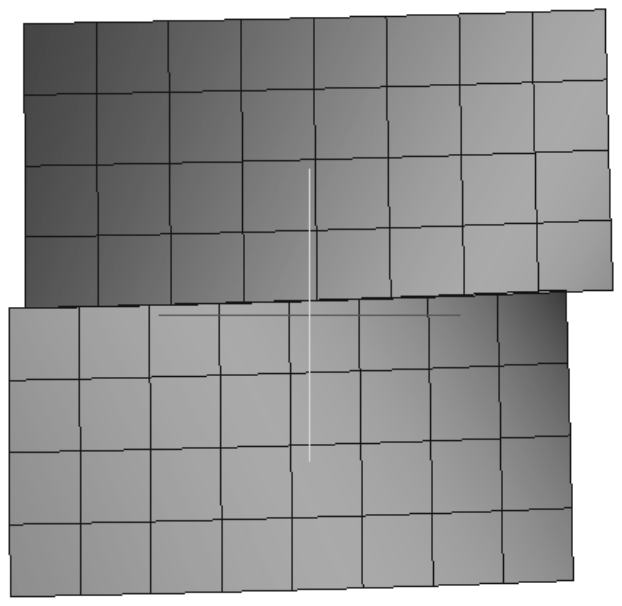}
\caption{Solution of the limit problem.} 
\label{fig:elastic03m8}
\end{figure}

Models used to simulate ground motion during an earthquake frequently
represent the sub--surface as a union of linearly elastic materials
separated by thin (fault) regions within which large deformations
(rupture) occur. Below we analyze limiting models which circumvent the
numerical difficulties encountered with direct simulation of these
models which arise when very fine meshes are required to resolve the
large deformations in the fault region. The fault region in reduced
models is represented as a surface and the rupture is realized as
discontinuities in certain components of the solution. Figures
\ref{fig:elastic3m844} and \ref{fig:elastic03m8} illustrates these
issues; the fine mesh in Figure \ref{fig:elastic3m844} is unnecessary
when the large shear across the fault is represented as the
discontinuity in the horizontal displacement shown in Figure
\ref{fig:elastic03m8}. Theorems \ref{thm:IepsToI} and \ref{thm:epsTo0}
justify this approach for a certain class of these models by
establishing that their solutions converge to the solution of a
reduced problem as the width of the fault tends to zero.

In the Appendix we show that a solution of the three dimensional
rupture model proposed in \cite{ZhAcWaBi15} can be found by considering a
cross section $\Omega = (-1,1) \times (-1,1)$ of a sub--surface region
containing a horizontal fault $S_\epsilon = (-1,1) \times
(-\epsilon/2, \epsilon/2)$ of width $\epsilon > 0$.
The balance of momentum takes the classical form
\begin{equation} \label{eqn:momentum}
\rho u_{tt} - \dv(T) = \rho f
\qquad \text{ on }
\qquad
(0,T) \times \Omega,
\end{equation}
where $u(t,x) \in \Re^2$ and $\rho(t,x) > 0$ represent the
displacement and density of the medium, $f(t,x) \in \Re^2$ the force
per unit mass, and $T$ is the (Cauchy) stress tensor. The sub--surface
strata is taken to be isotropic so that away from the fault the stress
takes the form
\begin{equation}
T = 2 \mu D(u) + \lambda \, \dv(u) I,
\qquad \text{ in } \Omega \setminus \Sbar_\epsilon,
\end{equation}
where
$D(u) = (1/2) \left(\nabla u + (\nabla u)^\top \right)$
and $\mu = \mu(x)$ and $\lambda = \lambda(x)$ are the shear
and bulk moduli. In the fault the stress is given by
\begin{equation} \label{eqn:T}
T = \begin{bmatrix}
      (2\mu+\lambda) u_{1x}   + \lambda u_{2y}&
      \mu \left(\epsilon (u_{2x} + u_{1y}) - \gamma  \vph\right) \vph \\
      \mu \left(\epsilon (u_{2x} + u_{1y}) - \gamma  \vph\right) &
      (2\mu+\lambda) u_{2y}   + \lambda u_{1x} \vph
    \end{bmatrix}
\qquad \text{ in } S_\epsilon,
\end{equation}
where $\gamma= \gamma(t,x)$ models the permanent deformation due to
damage, defects, and healing in the fault \cite{ZhAcWaBi15}, and
evolves according to
\begin{equation} \label{eqn:gamma}
(1/\beta) \gamma_t 
+ \etahat \gamma - T_{12} - \nu \gamma_{xx} = 0,
\qquad \text{ in } S_\epsilon.
\end{equation}
This model of rupture was inspired by the plasticity theories
developed in \cite{AcZh15,Ac10} where the coefficients $\beta$,
$\etahat$ and $\nu$ are typically nonlinear functions of $\gamma$ and
its derivatives.  To date there is no satisfactory mathematical theory
for these models of nonlinear elasticity \cite{AcTa11}, so for the
analysis below we assume that the coefficients in \eqnref{:gamma} are
specified which, for example, would be the case for one step of a
linearly implicit time stepping scheme for the fully nonlinear
problem. In this context we address the following problems.

{\em Strain Energy:} If $\gamma(x)$ is specified (or more generally
$\gamma_\epsilon \rightarrow \gamma$ sufficiently strongly) we verify
that the strain energies for the stationary problem
\eqnref{:momentum}--\eqnref{:T} converge to the limiting energy
$$
I(u) = \half \int_{\Omega \setminus S_0} 2\mu |D(u)|^2 + \lambda \, \dv(u)^2
+ \half \int_{S_0} \mu ([u_1]-\gamma)^2,
$$
where $[u_1(x)]$ denotes the jump in the first (horizontal) component
of $u$ across the line $S_0 = (-1,1) \times \{0\}$. The corresponding
Euler Lagrange operator is then
$$
-\dv(T) \quad \text{ on } \Omega\setminus S_0
\qquad \text{ with } \qquad
[Tn] = 0 \text{ and } T_{12} = \mu([u_1]-\gamma) 
\quad \text{ on } S_0,
$$
where $T = 2 \mu D(u) + \lambda \,\dv(u) I$ and $[Tn]$ is the jump of
the traction across $S_0$ with $n = (0,1)^\top$ and $T_{12}$ its
first component.

This shows that the scaling introduced in \eqnref{:T} is the
``mathematically interesting'' case for which a non--trivial limit
exits. With different scalings the equations for $u$ and $\gamma$
either decouple (the last term in the energy vanishes) or lock, $[u_1]
= \gamma$, in the limit.  The limiting energy for the coupled
stationary problem \eqnref{:momentum}--\eqnref{:gamma} is
\begin{equation} \label{eqn:coupled}
I(u,\gamma) 
= \half \int_{\Omega \setminus S_0} 2\mu |D(u)|^2 + \lambda \, \dv(u)^2 
+ \half \int_{S_0} \mu ([u_1]-\gamma)^2 
+ \left(\etahat \gamma^2 + \nu \gamma_x^2\right).
\end{equation}
We omit the proof since the proof is a routine extension of the
ideas used for the uncoupled problem. Examples of numerical solutions
to both the uncoupled and coupled problems are presented in Section 
\ref{sec:NumericalExamples}.

{\em Evolutionary Problem:} In Section \ref{sec:Evolution} we show
that solutions of the coupled system \eqnref{:momentum}--\eqnref{:gamma}
converge to a limit which satisfies the reduced system,
$$
\rho u_{tt} - \dv( T ) = \rho f 
\qquad \text{ on } \Omega \setminus S_0,
$$
with $T = 2 \mu D(u) + \lambda \,\dv(u) I$, and
\begin{equation*}
[Tn] = 0, \quad
T_{12} = \mu( [u_1] - \gamma ),
\text{ and }
(1/\beta) \gamma_t 
+ \etahat \gamma - T_{12} - \nu \gamma_{xx} = 0
\qquad \text{ on } S_0.
\end{equation*}

For definiteness we consider displacement boundary conditions $u(.,\pm
1)=0$ on the top and bottom of $\Omega$ and traction free boundary data
on the sides; $T(\pm 1, .) n = 0$, and $\nu \gamma_x(\pm 1,.) = 0$.
We omit analogous results for other boundary conditions which are
routine technical extensions of the proof techniques presented below.
The same energy and limiting problem are obtained with the
``engineering approximation'' utilized in \cite{ZhAcWaBi15} where the
shear stress $T_{12}(x,y)$ in the equation \eqnref{:gamma} is
approximated by its average $\Tbar_{12}(x)$ across $S_\epsilon$ so
that $\gamma$ depends only upon $x$,
$$
(1/\beta) \gamma_t 
+ \etahat \gamma -\Tbar_{12} - \nu \gamma_{xx}
= 0,
\quad \text{ where } \quad
\Tbar_{12}(t,x) 
= \frac{1}{\epsilon} \int_{-\epsilon/2}^{\epsilon/2} T_{12}(t,x,y) \, dy.
$$
The ideas presented below extend directly to the analysis of this
variation of the problem.

\subsection{Overview of Related Results}
The modeling and prediction of material failure is a notoriously
difficult problem with a long history. Broadly viewed, the onset of
failure is modeled either by crack formation (brittle failure) or
modifications of elasticity theory to admit plastic deformation
(ductile failure). Equations \eqnref{:momentum}--\eqnref{:gamma}
combine these approaches in the sense that a ``crack like'' fault
region is known a--priori; however, the material response in this zone
is inspired by models of dislocation motion in ductile materials.
This contrasts with ``cohesive zone'' models where friction laws
complement the usual jump condition, $[Tn] = 0$, and kinematic
compatibility condition, $[u_2] \geq 0$, (or $[u_2]=0$ on a fault deep
underground).  A detailed development of this model and comparision
with, and references to, alternative models may be found in
\cite[Section 8.4]{ZhAcWaBi15}.

A major challenge when modeling of brittle failure is predicting
nucleation and propagation paths of the cracks \cite{Horii1986brittle}. To date the
most mathematically complete approach to this issue is the variational
technique initiated by Ambrosio and Braides \cite{AmBr95}, which uses
spaces of bounded variation, and the corresponding extension to
evolution problems by Francfort and Marigo \cite{FrMa98} using the
concept of minimizing movements \cite{DeGiorgi93}. This has resulted in a
large body of mathematical work \cite{Ambrosio1995minimizing, DaTo02,FrLa03, Alberti2003calibration, Bourdin2008variational} and
numerical schemes \cite{Dolbow1999finite,Bourdin2000797, Areias2005analysis, Bourdin2007numerical, MeBoKh15}.  To date models based upon
these variational techniques do not include friction at the crack
sites or restrict the relative motion to prevent interpenetration.

Almost all models of ductile failure are based upon the premise that
the permanent macroscopic deformations of ductile materials result from the
motion of dislocations through a microscopic crystalline lattice
\cite{HiLo82}. However, passage from the microscopic description to a
macroscopic model with predictive capability requires a substantial
amount of phenomenological input \cite{Ac10,GuFrAn10}. Again the most
mathematically complete results for these models involve variational
techniques and quasi--static formulations to construct minimizing
movements \cite{Mi03,Dal2005quasistatic, Vi12}. The analogue of the friction laws in this
context is an assumption on the rate independence of the dissipation
which gives rise to energies with linear growth for which minimizers
have gradients of bounded deformation \cite{Mielke2004rate, Francfort2006existence,Dal2006quasistatic}.

%{\bf Jing: Look up more works on this: more of Mielke's work, people
%who have used his work, e.g. Theil; numerical work using Visintin or Mielke 
%ideas ...}

\subsection{Notation \& Function Spaces}

Standard notation is adopted for the Lebesgue spaces, $\Lp$, and the
Sobolev space $H^1(\Omega)$.  Solutions of evolution equations will be
viewed as functions from $[0,T]$ into these spaces, and we adopt the
usual notion, $\LtwoHone$, $\CHone$, etc.  to indicate the temporal
regularity. Strong and weak convergence in these spaces is denoted as
$u^\epsilon \rightarrow u$ and $u^\epsilon \weak u$ respectively.

Divergences of vector and matrix valued functions are denoted $\dv(u)
= u_{i,i}$ and $\dv(T)_i = T_{ij,j}$ respectively.  Here indices after
the comma represent partial derivatives and the summation convention
is used.  Gradients of vector valued quantities are interpreted as
matrices, $(\nabla u)_{ij} = u_{i,j}$, and the symmetric part of the
gradient is written as $D(u)$. Inner products are typically denoted as
pairings $(.,.)$ or, for clarity, the dot product of two vectors $v$,
$w \in \Re^d$ may be written as $v.w = v_i w_i$ and the Frobenious
inner product of two matrices $A$, $B \in \Re^{d \times d}$ as $A: B =
A_{ij} B_{ij}$.

The following notation is used to characterize the dependence upon
$\epsilon$ of the elastic and fault regions.

\begin{notation} \label{not:HoneEps}
  Let $\Omega = (-1,1)^2$ and $0 < \epsilon < 1/2$.
  \begin{enumerate} 
  \item The fault regions are denoted by $S_\epsilon = (-1,1) \times
    (-\epsilon/2, \epsilon/2)$ and $S_0 = (-1,1) \times \{0\}$ and
    their complements, the elastic regions, denoted as $
    \Omega_\epsilon = \Omega \setminus \Sbar_\epsilon$ and $\Omega_0 =
    \Omega \setminus S_0$.

  \item The sub--spaces of functions on the elastic region which vanish on
    the top and bottom boundaries are
    \begin{eqnarray*}
      U &=& 
      \{u \in H^1(\Omega) \sst u(x,\pm 1) = 0, \,\, -1 < x < 1\}, \\
      U_\epsilon &=& 
      \{u \in H^1(\Omega_\epsilon) \sst u(x,\pm 1) = 0, \,\, -1 < x < 1\}, \\
      U_0 &=& \{u \in H^1(\Omega_0) \sst u(x,\pm 1) = 0, \,\, -1 < x < 1\}. 
    \end{eqnarray*}

  \item The restriction $u \mapsto u|_{\Omega_\epsilon}$ is identified
    as an embedding of the spaces $\Hone \embed H^1(\Omega_\epsilon)$
    and $U \embed U_\epsilon$; similarly $\Hone \embed
    H^1(\Omega_0)$ and $U \embed U_0$.

  \item Below $\chi_A$ denotes the characteristic function of $A
    \subset \Omega$; $\chi_A(x) = 1$ if $x \in A$ and $\chi_A(x) = 0$ 
    otherwise.
  \end{enumerate}
\end{notation}

The following lemma quantifies the dependence upon $\epsilon$ of
embedding constants and properties of the function spaces for which
the energy is continuous and coercive. Here and below $C$ and $c$
denote constants which may vary from instance to instance but will
always be independent of $\epsilon$.

\begin{lemma} \label{lem:HoneEps}
  Denote the domains and spaces as in Notation \ref{not:HoneEps}, and if
  $u^\epsilon \in H^1(\Omega_\epsilon)$ and $u \in H^1(\Omega_0)$ denote
  by $[u^\epsilon]$ and $[u]$ the jump in their traces across the fault
  regions;
  $$ 
  [u^\epsilon] = u^\epsilon(.,\epsilon/2) - u^\epsilon(.,-\epsilon/2)
  \qquad \text{ and } \qquad
  [u] = u(.,0^+) - u(.,0^-).
  $$
  \begin{enumerate}
  \item The constant in Korn's inequality on $U_\epsilon$ is
    independent of $\epsilon$.
    
  \item The following Poincare inequality holds for functions in $U$,
    $$
    (1/2) \norm{u}_{L^2(S_\epsilon)}
    \leq \left(\epsilon^2 \norm{u_y}_{L^2(S_\epsilon)}^2
      +  (\epsilon/2) \norm{u}_{H^1(\Omega_\epsilon)}^2 \vph \right)^{1/2}.
    $$
    
  \item If $u^\epsilon \in \Hone$ and
    $$
    u^\epsilon \weak u, \qquad 
    \chi_{\Omega_\epsilon} u^\epsilon_x \weak g_0 \qquad
    u^\epsilon_y \weak g_1, \qquad
    \text{ in } \Ltwo,
    $$
    then $u \in \Hone$ and $\nabla u = (g_0, g_1)^\top$.
    
  \item If $u^\epsilon \in \Hone$ and
    $$
    u^\epsilon \weak u, \qquad 
    u^\epsilon_x \weak g_0, \qquad
    \chi_{\Omega_\epsilon} u^\epsilon_y \weak g_1 \qquad
    \text{ in } \Ltwo,
    $$
    then $u \in H^1(\Omega_0)$ and $\nabla u = (g_0, g_1)^\top$. In addition
    $[u^\epsilon] \rightarrow [u]$ in $L^2(-1,1)$.
    
  \item Let $\phi_\epsilon:\Omega_\epsilon \rightarrow \Omega_0$ be the
    mapping
    $$
    \phi_\epsilon(x,y) = 
    \left(x, \tfrac{y-\epsilon/2}{1-\epsilon/2} \right),
    \quad \epsilon/2 < y < 1, 
    \quad \text{ and } \quad
    \phi_\epsilon(x,y) = 
    \left(x, \tfrac{y+\epsilon/2}{1-\epsilon/2} \right),
    \quad -1 < y < -\epsilon/2.
    $$
    Then the linear functions $E_\epsilon: H^1(\Omega_\epsilon)
    \rightarrow H^1(\Omega_0)$ given by $E_\epsilon(u^\epsilon) =
    u^\epsilon \circ \phi_\epsilon^{-1}$ are isomorphisms and their
    norms and the norms of their inverses converge to one as $\epsilon
    \rightarrow 0$. The restriction of $E_\epsilon$ to $U_\epsilon$ is
    an isomorphism onto $U_0$.
    
  \item \label{it:extension}
    If $u \in H^1(\Omega_0)$ then there exists $u^\epsilon \in \Hone$
    such that
    \begin{enumerate}
    \item $\norm{E_\epsilon(u^\epsilon) - u}_{H^1(\Omega_0)} \rightarrow 0$.
      
    \item $\norm{u^\epsilon_x}_{L^2(S_\epsilon)} \rightarrow 0$.
      
    \item $u^\epsilon_{1y}(x,.)$ is independent of $y$ in $S_\epsilon$ and
      $\int_{-\epsilon/2}^{\epsilon/2} u^\epsilon_{1y}(.,y) \, dy = 
      [u^\epsilon] \rightarrow [u]$ in $L^2(-1,1)$.
    \end{enumerate}
    In addition, if $u \in U_0$ then $u^\epsilon \in U$.
  \end{enumerate}
\end{lemma}

\begin{figure} 
\begin{center}
\includegraphics[height=1.25in]{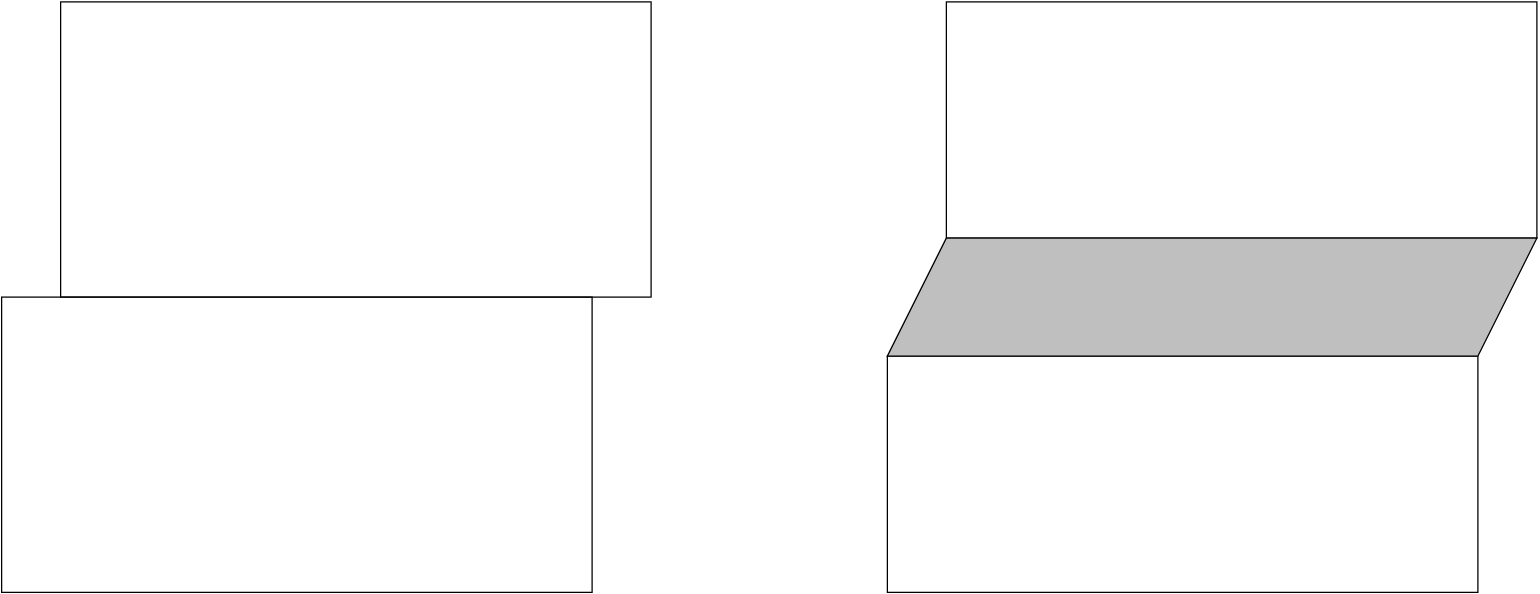}
\caption{Approximation of horizontal displacement $u \in U_0$ by
  a function $u^\epsilon \in U \subset \Hone$} \label{fig:U0toU}
\end{center}
\end{figure}

The proof of this lemma involves standard arguments \cite{Sh77} so is
omitted. The only subtlety appears in the construction of the function
$u^\epsilon$ in item \ref{it:extension}. The idea is illustrated in
Figure \ref{fig:U0toU}; given $u \in U_0$ extend $u \circ
\phi^{-1}_\epsilon$ to $S_\epsilon$ by linear interpolation.  However,
to control the $x$-derivative in $S_\epsilon$ we first mollify $u$ on
each subdomain $\Omega^\pm$ with parameter $\sqrt{\epsilon}$.

\section{Gamma Convergence of the Stationary Operator} \label{sec:Gamma} 

Letting $\bbC:\Re^{2 \times 2} \rightarrow \Re^{2 \times 2}_{sym}$
denote the classical isotropic elasticity tensor with shear and bulk
moduli $\mu$ and $\lambda$, the associated strain energy function will
be denoted as
$$
2 W(A) \equiv \bbC(A):A = 2\mu (A_{11}^2 + A_{22}^2) 
+ \lambda (A_{11} + A_{22})^2 
+ \mu (A_{12} + A_{21})^2.
$$
Define $I_\epsilon: \Hone^2 \rightarrow \Re$ to be energy,
\begin{eqnarray}
  I_\epsilon(u) 
  &=& \int_{\Omega_\epsilon}
  W\left(\nabla u \vph\right)
  + \int_{S_\epsilon} W\left(
    \begin{bmatrix} 
      u_{1x} & \sqrt{\epsilon} u_{1y} - \gamma/\sqrt{\epsilon} \\
      \sqrt{\epsilon} u_{2x} & u_{2y}
    \end{bmatrix} \right)  \label{eqn:Ieps} \\ [1ex]
  &=& \int_\Omega W \left(
    \begin{bmatrix} 
      u_{1x} & u_{1y} \chi_{\Omega_\epsilon} \\
      u_{2x} \chi_{\Omega_\epsilon} & u_{2y}
    \end{bmatrix} \right)
  + \int_{S_\epsilon} W\left(
    \begin{bmatrix} 
      0 & \sqrt{\epsilon} u_{1y} - \gamma/\sqrt{\epsilon} \\
      \sqrt{\epsilon} u_{2x} & 0
    \end{bmatrix} \right) 
  \nonumber
\end{eqnarray}

In this section we set up the following theorem which establishes
convergence of the energies in the following sense \cite{DaM93}.
\begin{itemize}
\item (lim--inf inequality) If $\{u_\epsilon\}_{\epsilon > 0} \subset
  U \times U$ with $\{I(u_\epsilon)\}_{\epsilon > 0} \subset \Re$
  bounded then there exists $u \in U_0 \times U$ and a sub--sequence
  for which $u^\epsilon \weak u$ in $H^1_{loc}(\Omega_0) \times
  \Hone$ and $I(u) \leq \liminf_{\epsilon \rightarrow 0}
  I_\epsilon(u^\epsilon)$.

\item (lim--sup inequality) For each $u \in U_0 \times U$ there exists
  a sequence $\{u_\epsilon\}_{\epsilon > 0} \subset U \times U$ such
  that $u^\epsilon \weak u$ in $H^1_{loc}(\Omega_0) \times \Hone$ and
  $I(u) \geq \limsup_{\epsilon \rightarrow 0} I_\epsilon(u^\epsilon)$.
\end{itemize}

\begin{theorem} \label{thm:IepsToI} Denote the domains and spaces as
  in Notation \ref{not:HoneEps} and let $I_\epsilon:U \times U
  \rightarrow \Re$ be as in equation \eqnref{:Ieps} with $\gamma \in
  L^2(-1,1)$ fixed. Assume that the shear and bulk moduli are bounded
  above, there exists $c_0 > 0$ such that $\mu \geq c_0$ and $\mu +
  \lambda \geq c_0$, and that there exists $\epsilon_0 > 0$ such that
  the shear modulus $\mu$ is independent of $y$ on $S_\epsilon$ for
  $\epsilon < \epsilon_0$. Then $I_\epsilon
  \stackrel{\Gamma}{\rightarrow} I$ where $I:U_0 \times U \rightarrow
  \Re$ is given by
  $$
  I(u) = \int_{\Omega_0} W(\nabla u) 
  + \half \int_{-1}^1 \mu ([u_1] - \gamma)^2,
  $$ 
  for which the strong form of the  Euler Lagrange operator is
  $$
  -\dv(\bbC(\nabla u)) \text{ on } \Omega_0, 
  \quad \text{ with } \quad
  [\bbC(\nabla u)].n = 0
  \text{ and }
  \bbC(\nabla u)_{12} = \mu([u_1]-\gamma) 
  \text{ on } S_0.
  $$
  Here $[.]$ denotes the jump across the fault line $y=0$.
\end{theorem}

The following lemma quantifies the coercivity properties of the
energies $I_\epsilon$ and the corresponding bounds required for the
proof of Theorem \ref{thm:IepsToI}. In this lemma we use the property
that in two dimensions the assumptions on the Lame parameters
guarantee $W(A) \geq 2 c_0 |A_{sym}|^2$ where $A_{sym}$ denotes the
symmetric part of $A \in \Re^{2 \times 2}$.

\begin{lemma} \label{lem:Ibound}
  Denote the domains and spaces as in Notation \ref{not:HoneEps} and
  let $I_\epsilon:U \times U \rightarrow \Re$ be as in equation
  \eqnref{:Ieps} with $\gamma \in L^2(-1,1)$ fixed. Assume that the
  shear and bulk moduli bounded above, $\mu \geq c_0 > 0$ and
  $\mu + \lambda \geq c_0 > 0$. Then
  $$
  \ltwo{u_{1x}}^2
  + \norm{u_{1y}}_{L^2(\Omega_\epsilon)}^2
  + \norm{u_{2x}}^2_{L^2(\Omega_\epsilon)} 
  + \ltwo{u_{2y}}^2
  \leq C I_\epsilon(u),
  $$
  and
  $$
  \norm{\sqrt{\epsilon} u_{2x}}_{L^2(\Omega_0)}^2
  + \norm{\sqrt{\epsilon} u_{1y}}_{L^2(\Omega_0)}^2
  \leq C \left(I_\epsilon(u) + \norm{\gamma}_{L^2(-1,1)}^2 \right).
  $$
  In particular,
  $$
  \ltwo{u_1}^2 + \ltwo{u_2}^2
  \leq C \left(I_\epsilon(u) + \norm{\gamma}_{L^2(-1,1)}^2 \right).
  $$
\end{lemma}

\begin{proof}
  It is immediate that
  $$
  \ltwo{u_{1x}}^2 
  + \norm{u_{1y} + u_{2x}}_{L^2(\Omega_\epsilon)}^2
  + \ltwo{u_{2y}}^2
  + \norm{\sqrt{\epsilon}(u_{2x} + u_{1y})
    - \gamma/\sqrt{\epsilon}}_{L^2(S_\epsilon)}^2
  \leq C I_\epsilon(u),
  $$
  and Korn's inequality on $\Omega_\epsilon$ shows
  $$
  \norm{u_{2x}}_{L^2(\Omega_\epsilon)}^2
  + \norm{u_{1y}}_{L^2(\Omega_\epsilon)}^2
  \leq C \left( \norm{u_{1y} + u_{2x}}_{L^2(\Omega_\epsilon)}^2
  + \norm{u_{1x}}_{L^2(\Omega_\epsilon)}^2
  + \norm{u_{2y}}_{L^2(\Omega_\epsilon)}^2 \right).
  $$
  Next, use the triangle inequality and the identity
  $\norm{\gamma/\sqrt{\epsilon}}_{L^2(S_\epsilon)} =
  \norm{\gamma}_{L^2(-1,1)}$ to obtain
  $$
  \norm{\sqrt{\epsilon} (u_{2x} + u_{1y})}_{L^2(S_\epsilon)} 
  \leq
  \norm{\sqrt{\epsilon} (u_{2x} + u_{1y})
    - \gamma/\sqrt{\epsilon}}_{L^2(S_\epsilon)}
  + \norm{\gamma}_{L^2(-1,1)}, 
  $$
  Korn's inequality for the vector field 
  $\ut = \sqrt{\epsilon} (u_1, u_2)$
  on $\Omega_0$ shows
  \begin{multline*}
    \norm{\sqrt{\epsilon} u_{2x}}_{L^2(\Omega_0)}^2
    + \norm{\sqrt{\epsilon} u_{1y}}_{L^2(\Omega_0)}^2 \\
    \leq C \left(
      \norm{\sqrt{\epsilon}(u_{2x} +  u_{1y})}_{L^2(\Omega_0)}^2
      + \epsilon \big(\norm{u_{1x}}_{L^2(\Omega_0)}^2
      + \norm{u_{2y}}_{L^2(\Omega_0)}^2 \big) \right).
  \end{multline*}
\end{proof}

\begin{proof} (of Theorem \ref{thm:IepsToI})
{\em Lim--Inf Inequality:}
Let $\{u^\epsilon\}_{\epsilon > 0} \subset \Hone^2$ and
suppose $I_\epsilon(u^\epsilon)$ is bounded. Lemma \ref{lem:Ibound} 
then shows that the functions
$$
u^\epsilon_1, \,\, u^\epsilon_{1x}, \,\, u^\epsilon_{1y}  \chi_{\Omega_\epsilon},
\qquad \text{ and } \qquad
u^\epsilon_2, \,\, u^\epsilon_{2x} \chi_{\Omega_\epsilon}, \,\, u^\epsilon_{2y}, 
$$
are all bounded in $\Ltwo$. Upon passing to a subsequence we may then
assume each of them converges weakly in $\Ltwo$ to a limit $u =
(u_1,u_2) \in \Ltwo^2$ and from Lemma \ref{lem:HoneEps} conclude $u
\in U_0 \times U$; in particular,
$$
\begin{bmatrix} 
  u^\epsilon_{1x} & u^\epsilon_{1y} \chi_{\Omega_\epsilon} \\
  u^\epsilon_{2x} \chi_{\Omega_\epsilon} & u^\epsilon_{2y} \end{bmatrix} 
\weak \begin{bmatrix} 
  u_{1x} & u_{1y} \\
  u_{2x} & u_{2y} \end{bmatrix}
\qquad \text{ in } \Ltwo^{2 \times 2},
$$
Since $W$ is convex and continuous it is weakly lower
semi--continuous; in particular, the limit of the first term in
equation \eqnref{:Ieps} is bounded as
$$
\int_{\Omega_0} W\left( \nabla u \right)
\leq \liminf_\epsilon
\int_{\Omega} W\left( 
  \begin{bmatrix} 
    u^\epsilon_{1x} & u^\epsilon_{1y} \chi_{\Omega_\epsilon} \\
    u^\epsilon_{2x} \chi_{\Omega_\epsilon} & u^\epsilon_{2y} 
  \end{bmatrix} \right).
$$
To compute the limit of the second term in equation \eqnref{:Ieps},
use Jensen's
inequality and the quadratic homogeneity of $W(.)$ to obtain
\begin{eqnarray*}
  \int_{S_\epsilon} W\left( 
    \begin{bmatrix} 
      0 & \sqrt{\epsilon} u^\epsilon_{1y} - \gamma/\sqrt{\epsilon} \\
      \sqrt{\epsilon} u^\epsilon_{2x} & 0
    \end{bmatrix} \right)
  &\geq& 
  \int_{-1}^1 \epsilon W\left( 
    \frac{1}{\epsilon} \int_{-\epsilon/2}^{\epsilon/2}  
    \begin{bmatrix} 
      0 & \sqrt{\epsilon} u^\epsilon_{1y} - \gamma/\sqrt{\epsilon} \\
      \sqrt{\epsilon} u^\epsilon_{2x} & 0
    \end{bmatrix}
    \, dy \right) \, dx \\
  &=& 
  \int_{-1}^1  W\left(
    \begin{bmatrix}
      0 & [u^\epsilon_1] - \gamma  \\
      \int_{-\epsilon/2}^{\epsilon/2}  u^\epsilon_{2x} \, dy & 0
    \end{bmatrix} \right) \, dx,
\end{eqnarray*}
where $[u^\epsilon_1](x) = u_1(x,\epsilon/2) - u_1(x, -\epsilon/2)$.
Lemma \ref{lem:HoneEps} shows $[u^\epsilon_1] \rightarrow [u]$ 
in $L^2(-1,1)$, so the lim--inf inequality will follow upon showing that
$\int_{-\epsilon/2}^{\epsilon/2} u^\epsilon_{2x} \, dy \weak 0$
in $L^2(-1,1)$. To verify this, first use
the Cauchy Schwarz inequality and Lemma \ref{lem:Ibound} to bound
this term in $L^2(-1,1)$,
$$
\int_{-1}^1 \left(
  \int_{-\epsilon/2}^{\epsilon/2} u^\epsilon_{2x} \, dy \right)^2 \, dx
% \leq \int_{-1}^1 \int_{-\epsilon/2}^{\epsilon/2}
% (\sqrt{\epsilon} u^\epsilon_{2x})^2 \, dy \, dx
\leq \norm{\sqrt{\epsilon} u_{2x}^\epsilon}_{L^2(S_\epsilon)}^2
\leq C \left( I_\epsilon(u^\epsilon_1) + \norm{\gamma}_{L^2(-1,1)}^2 \right)^{1/2}.
$$
To show that this term converges weakly to zero let $\phi \in
C_0^\infty(-1,1)$ and compute
$$
\left| \int_{-1}^1 \int_{-\epsilon/2}^{\epsilon/2} 
  u^\epsilon_{2x} \, dy 
  \, \phi \, dx \right|
= \left| \int_{-1}^1 \int_{-\epsilon/2}^{\epsilon/2} 
  u^\epsilon_2 
  \, \phi' \, dy \, dx \right|
\leq
\norm{u^\epsilon_2}_{L^2(S_\epsilon)} \sqrt{\epsilon} \norm{\phi'}_{L^2(-1,1)} .
$$
The sharp Poincare inequality in Lemma \ref{lem:HoneEps} shows
$\norm{u^\epsilon_2}_{L^2(S_\epsilon)} \leq C \sqrt{\epsilon}$ so
the right hand side of the above vanishes from which it
follows that $ \int_{-\epsilon/2}^{\epsilon/2} u^\epsilon_{2x} \, dy
\weak 0$ in $L^2(-1,1)$.

{\em Lim--Sup Inequality:}
To construct a recovery sequence for $u \in H^1(\Omega_0) \times \Hone$ 
select $u^\epsilon = (u_1^\epsilon, u_2)$ where $u^\epsilon_1$
is the lifting of $u_1$ to $\Hone$ guaranteed by item \ref{it:extension}
of Lemma \ref{lem:HoneEps}. 

Lemma \ref{lem:HoneEps} shows $u^\epsilon \circ \phi_\epsilon^{-1}
\rightarrow u$ in $H^1(\Omega_0)$ where $\phi_\epsilon:\Omega_\epsilon
\rightarrow \Omega_0$ is the piecewise affine diffeomorphism in the
Lemma.  Since the mapping $u \mapsto W(\nabla u)$ is continuous on
$H^1(\Omega_0)^2$ it follows that the energy in the bulk converges,
$$
\int_{\Omega_\epsilon} W(\nabla u^\epsilon) 
=\int_{\Omega_0} W\left( \nabla (u^\epsilon \circ \phi_\epsilon^{-1}) \vph \right)
(1 - \epsilon/2) 
\rightarrow \int_{\Omega_0} W(\nabla u).
$$
The energy in the fault regions $S_\epsilon$ takes the form
$$
\int_{S_\epsilon} W \left(
  \begin{bmatrix} 
    u^\epsilon_{1x} 
    & \sqrt{\epsilon} u^\epsilon_{1y} - \gamma/\sqrt{\epsilon} \\
    \sqrt{\epsilon} u_{2x} 
    & u_{2y}
  \end{bmatrix} \right)
= \int_{S_\epsilon} W \left(
  \begin{bmatrix} 
    u^\epsilon_{1x} 
    & ([u^\epsilon_1] - \gamma)/ \sqrt{\epsilon} \\
    \sqrt{\epsilon} u_{2x} 
    & u_{2y}
  \end{bmatrix} \right)
$$
Since $u_{2x}$, $u_{2y} \in \Ltwo$ are independent of $\epsilon$ and
$|S_\epsilon| \rightarrow 0$ it is immediate that
$\norm{u_{2x}}_{L^2(S_\epsilon)}$ and
$\norm{u_{2y}}_{L^2(S_\epsilon)}$ both converge to zero,
and from Lemma \ref{lem:HoneEps} it follows that
$\norm{u^\epsilon_{1x}}_{L^2(S_\epsilon)}$ also converges to zero.
Also, $[u^\epsilon_1] - \gamma$ is independent of $y$ and
$[u^\epsilon_1] \rightarrow [u_1]$ in $L^2(-1,1)$ so
$$
\frac{1}{\sqrt{\epsilon}} \norm{[u^\epsilon_1] - \gamma}_{L^2(S_\epsilon)}
=  \norm{[u^\epsilon_1] - \gamma}_{L^2(-1,1)}  
\rightarrow \norm{[u_1] - \gamma}_{L^2(-1,1)}.
$$
Since $W:\Re^{2 \times 2} \rightarrow \Re$ is continuous, non--negative, 
and has quadratic growth it follows that
$$
\int_{S_\epsilon} W \left(
  \begin{bmatrix} 
    u^\epsilon_{1x} 
    & \sqrt{\epsilon} u^\epsilon_{1y} - \gamma/\sqrt{\epsilon} \\
    \sqrt{\epsilon} u_{2x} 
    & u_{2y}
  \end{bmatrix} \right)
\rightarrow \int_{S_0} W \left(
  \begin{bmatrix} 
    0 & ([u_1] - \gamma) \\ 0 & 0
  \end{bmatrix} \right),
$$
and $\{u^\epsilon\}_{\epsilon>0} \subset \Hone^2$ is a recovery sequence.
\end{proof}

\section{Evolution Equations} \label{sec:Evolution} In this section we
show that as $\epsilon \rightarrow 0$ solutions of equations
\eqnref{:momentum}--\eqnref{:gamma} converge to the solution of a
limiting problem with the spatial Euler Lagrange operator
corresponding to the gamma limit obtained in the previous section.

Solutions of equations \eqnref{:momentum}--\eqnref{:gamma} satisfy
$(u(t), \gamma(t)) \in U^2 \times G_\epsilon$ and
\begin{gather}
\int_\Omega (\rho u_{tt}, \uhat) 
+ \left( \bbC_\epsilon(D(u)), D(\uhat) \right)
+ \int_{S_\epsilon} \mu \left(
  \epsilon (u_{2x} + u_{1y}) 
  - \gamma, \uhat_{2x}+\uhat_{1y} \vph\right)
= \int_\Omega (\rho f,\uhat), \label{eqn:uEpsT} \\
\frac{1}{\epsilon} \int_{S_\epsilon} 
(1/\beta) (\gamma_t, \gammahat) 
+ \ell\left(\gamma,\gammahat \right)
- \mu \left(\epsilon(u_{2x} 
  + u_{1y})-\gamma, \gammahat \right) = 0, \label{eqn:gEpsT}
\end{gather}
for all $(\uhat,\gammahat) \in U^2 \times G_\epsilon$ where
$$
U = \{u \in \Hone \sst u(.,\pm 1) = 0 \}
\qquad \text{ and } \qquad
G_\epsilon = \{\gamma \in L^2(S_\epsilon) \sst \gamma_x \in L^2(S_\epsilon)\}.
$$
In this weak statement $D(u) = (1/2)(\nabla u + (\nabla u)^\top)$ is
the symmetric part of the displacement gradient and
$$
\bbC_\epsilon(D) = \bbC\left( 
\begin{bmatrix} D_{11} & D_{12} \chi_{\Omega_\epsilon} \\
D_{21} \chi_{\Omega_\epsilon} & D_{22}
\end{bmatrix} \right),
\qquad
\ell(\gamma,\gammahat) = \nu \gamma_x \gammahat_x + \etahat\gamma \gammahat,
$$
where $\bbC(D) = 2 \mu D + \lambda \,tr(D) I$ is the isotropic
elasticity tensor.

Solutions of the sharp interface problem satisfy $(u(t), \gamma(t)) \in (U_0 \times U)
\times G$ and
\begin{gather}
\int_{\Omega_0} (\rho u_{tt}, \uhat) + \left( \bbC(D(u)), D(\uhat) \right)
+ \int_{-1}^1 \mu \left([u_1] - \gamma, [\uhat_1] \vph\right)
= \int_\Omega (\rho f,\uhat), \label{eqn:u0T} \\
\int_{-1}^1 (1/\beta) (\gamma_t, \gammahat) 
+ \ell\left(\gamma,\gammahat \right)
- \mu \left([u_1]-\gamma, \gammahat \right) = 0, \label{eqn:g0T}
\end{gather}
for all $(\uhat,\gammahat) \in (U_0 \times U) \times G$ where
$$
U_0 = \{u \in H^1(\Omega_0) \sst u(.,\pm 1) = 0\}
$$
$$
U = \{u \in \Hone \sst u(.,\pm 1) = 0\}
\qquad \text{ and } \qquad
G = H^1(-1,1).
$$

In this section we prove the following theorem which establishes
convergence of solutions of equations \eqnref{:uEpsT}--\eqnref{:gEpsT}
to solutions of \eqnref{:u0T}--\eqnref{:g0T}.

\begin{theorem} \label{thm:epsTo0} Denote the domains and spaces as in
  Notation \ref{not:HoneEps} and assume that the coefficients in
  equations \eqnref{:uEpsT}--\eqnref{:gEpsT} are independent of time
  and there exist constants $C$, $c$ such that
  $$
  0 < c \leq \rho(x), \, \mu(x), \, \beta(x), \, \nu(x), \, 
  \mu(x) + \lambda(x) \leq C, 
  \quad \text{ and } \quad
  0 < \etahat(x) < C,
  $$
  and that there exists $\epsilon_0 > 0$ such that the shear modulus
  $\mu$ is independent of $y$ on $S_\epsilon$ for $\epsilon <
  \epsilon_0$. 

  Fix $f \in L^1[0,T;\Ltwo]$ and initial data $u_t(0) \in \Ltwo$, $\gamma(0)
  \in H^1(-1,1)$ and $u(0) \in U_0 \times U$ for the sharp interface
  problem and let the initial values for equations
  \eqnref{:uEpsT}--\eqnref{:gEpsT} be
  $$
  u^\epsilon_t(0) = u_t(0), \qquad 
  u^\epsilon_2(0)=u_2(0), \qquad
  \gamma^\epsilon(0) = \gamma(0),
  $$
  and $u^\epsilon_1(0) = u_1$ if $u_1(0) \in U$; otherwise,
  select $\{u_1^\epsilon(0)\}_{\epsilon>0} \subset U$ such that
  $$
  \norm{u_1^\epsilon(0) - u_1(0)}_{H^1(\Omega_\epsilon)} \rightarrow 0
  \quad \text{ and } \quad
  \norm{u^\epsilon_{1x}(0)}_{L^2(S_\epsilon)} 
  + \sqrt{\epsilon} \norm{u^\epsilon_{1y}(0)}_{L^2(S_\epsilon)} 
  \leq C \norm{u_1(0)}_{H^1(\Omega_0)}.
  $$
  Let $(u^\epsilon, \gamma^\epsilon)$ denote the solution of
  \eqnref{:uEpsT}--\eqnref{:gEpsT} with this data and let
  $\gammabar^\epsilon(t,x) = (1/\epsilon)
  \int_{-\epsilon/2}^{\epsilon/2} \gamma^\epsilon(t,x,y) \, dy$. Then
  $\{(u^\epsilon, \gammabar^\epsilon)\}_{\epsilon>0}$ converges weakly
  in $H^1[0,T; \Ltwo^2 \times L^2(-1,1)]$ and strongly in $L^2[0,T;
  \Ltwo^2 \times L^2(-1,1)]$ to a limit 
  $$
  (u,\gamma) \in \calU \equiv H^1[0,T;\Ltwo^2 \times L^2(-1,1)] 
  \cap L^2[0,T; (U_0 \times U) \times H^1(-1,1)]
  $$ 
  with initial data $(u(0), \gamma(0))$ which satisfies
  \begin{gather*}
    \int_0^T\!\int_{\Omega_0} -(\rho u_t, \uhat_t) 
    + \left( \bbC(D(u)), D(\uhat) \right)
    + \int_0^T\!\int_{-1}^1 \mu \left( [u_1] 
      - \gamma,  [\uhat_1] \right) 
    = \int_\Omega (\rho u_t(0), \uhat(0)) 
    +\int_0^T\!\int_\Omega (\rho f,\uhat), \\
    \int_0^T\!\int_{-1}^1 
    (1/\beta) (\gamma_t, \gammahat) 
    + \ell\left(\gamma,\gammahat \right)
    - \mu \left([u_1]-\gamma, \gammahat \right) 
    = 0,
  \end{gather*}
  for all $ (\uhat,\gammahat) \in \calU$ with $\uhat(T) = 0$.
\end{theorem}

\subsection{Existence of Solutions and Bounds}

Equations \eqnref{:uEpsT}--\eqnref{:gEpsT} and
\eqnref{:u0T}--\eqnref{:g0T} both have the structure of a degenerate
wave equation on a product spaces taking the form $(u(t),\gamma(t))
\in \bbU \times \bbG$,
\begin{equation} \label{eqn:wave0}
C(u,\gamma)_{tt} + B (u,\gamma)_t + A(u,\gamma) = (\rho f,0),
\end{equation}
with
\begin{equation} \label{eqn:CB}
C(u,\gamma) = (\rho u, 0), 
\qquad \text{ and } \qquad
B(u,\gamma) = (0,\gamma/\beta),
\end{equation}
(the later scaled by $1/\epsilon$ for the $\epsilon$ equation)
and $A:\bbU \times \bbG \rightarrow \bbU' \times \bbG'$ is the Riesz map for the
space $\bbU \times \bbG$.  For the limit problem
$$
A(u,\gamma)(u,\gamma) = \norm{(u,\gamma)}_0^2
= \int_{\Omega_0} \bbC(D(u)):D(u)
+ \int_{-1}^1 \ell\left(\gamma,\gamma) + \mu ([u_1] - \gamma \vph\right)^2,
$$
and for the $\epsilon$ equation
$A(u,\gamma)(u,\gamma) = \norm{(u,\gamma)}^2_\epsilon$ with
$$
\norm{(u,\gamma)}^2_\epsilon
= \int_\Omega \bbC_\epsilon(D(u)):D(u)
+ \int_{S_\epsilon} (1/\epsilon) \ell(\gamma,\gamma)         
+ \mu \left(\sqrt{\epsilon}(u_{2x} + u_{1y}) 
  - \gamma/\sqrt{\epsilon} \vph\right)^2.
$$
The hypotheses on the initial data in Theorem \ref{thm:epsTo0}
guarantee $\norm{(u^\epsilon(0),\gamma^\epsilon)(0)}_\epsilon \rightarrow
\norm{(u(0),\gamma(0))}_0$.

The following theorem from \cite[Corollary VI.4.2]{Sh77} establishes
existence of (strong) solutions to equations which take the form shown
in \eqnref{:wave0}. In the statement of this theorem $\calL(V,V')$
denotes the continuous linear operators and $B \in \calL(V,V')$ is
monotone if $Bv(v) \geq 0$ for all $v \in V$.

\begin{theorem} \label{thm:wave} Let $A$ be the Riesz map of the
  Hilbert space $V$ and let $W$ be the semi--normed space obtained
  from the symmetric and monotone $C \in \calL(V,V')$. Let $D(B)
  \subset V$ be the domain of a linear monotone operator $B:D(B)
  \rightarrow V'$. Assume that $B+C$ is strictly monotone and
  $A+B+C:D(B) \rightarrow V'$ is surjective. Then for every $f \in
  C^1[0,\infty,W')$ and every pair $v_0 \in V$ and $v_1 \in D(B)$ with
  $Av_0 + Bv_1 \in W'$, there exists a unique 
  $$
  v \in C[0,\infty,V) \cap C^1(0,\infty,V) \cap
  C^1[0,\infty,W) \cap C^2(0,\infty,W),
  $$ 
  with $v(0)=v_0$, $Cv'(0)=Cv_1$ and for each $t > 0$, $v' \in
  D(B)$, $Av(t) + Bv'(t) \in W'$ and
  \begin{equation} \label{eqn:wave}
    (Cv'(t))' + B v'(t) + Av(t) = f(t).
  \end{equation}
\end{theorem}

When $V = U^2 \times G_\epsilon$ or $(U_0 \times U) \times G$ with
operators as in equation \eqnref{:CB} the state space is $W =
L^2(\Omega)$ with weight $\rho$ and $D(B) = V$ is the whole space.
Then $(B+C)(u,\gamma) = (\rho u, \gamma / \beta)$ is strictly
monotone, and $C+B+A:V \rightarrow V'$ is the sum of the Riesz map
with a monotone map, so is surjective.

The existence of strong solutions guaranteed by Theorem \ref{thm:wave} was
obtained upon writing equation \eqnref{:wave} as a first order system,
$\calB(v,v')' + \calA(v,v') = \ft$ with
$$
\calB = \begin{bmatrix} A & 0 \\ 0 & C \end{bmatrix} 
\qquad 
\calA = \begin{bmatrix} 0 & -A \\ A & B \end{bmatrix}
\qquad \text{ and } \qquad
\ft = \begin{pmatrix} 0 \\ f \end{pmatrix}.
$$
Classical semi--group theory then provides necessary and
sufficient conditions upon the data for the existence of strong
solutions. An alternative to the semi--group approach is to use
\cite[Proposition III.3.3]{Sh97} which establishes existence of weaker
solutions for a broader class of data and problems with time dependent
coefficients. Weak solutions exist when $f \in L^1[0,T;W']$ and
satisfy
$$
|v'|_{L^\infty[0,T;W]}^2 + \norm{v}_{L^\infty[0,T;V]}^2 + \int_0^T Bv'(v') 
\leq C \left( |v'(0)|_W^2 + \norm{v(0)}_V^2 
+ \norm{f}_{L^1[0,T;W']}^2 \right).
$$
The following corollary summarizes bounds available for solutions 
of \eqnref{:uEpsT}--\eqnref{:gEpsT} that results from this theory
and the Korn and sharp Poincare inequalities stated in Lemma \ref{lem:HoneEps}.

\begin{corollary} \label{cor:bounds}
  Under the hypotheses of Theorem \ref{thm:epsTo0} there exists
  a constant $C > 0$ independent of $\epsilon$ for which solutions
  $(u^\epsilon, \gamma^\epsilon)$ of \eqnref{:uEpsT}--\eqnref{:gEpsT}
  satisfy 
  \begin{multline*}
%    \norm{u^\epsilon_{tt}}_{L^2[0,T;H^{-1}(\Omega_\epsilon)]} +
    \linfltwo{u^\epsilon_t} 
    + \norm{u^\epsilon}_{L^\infty[0,T;H^1(\Omega_\epsilon)]} \\
    \hspace*{0.25in}  
    + (1/\sqrt{\epsilon}) \norm{\gamma^\epsilon_t}_{L^2[0,T;L^2(S_\epsilon)]} 
    + (1/\sqrt{\epsilon}) \norm{\gamma^\epsilon}_{L^\infty[0,T;L^2(S_\epsilon)]} 
    + (1/\sqrt{\epsilon}) \norm{\gamma^\epsilon_x}_{L^\infty[0,T;L^2(S_\epsilon)]} 
    \hfill \\
    + \norm{u_{1x}}_{L^\infty[0,T;L^2(S_\epsilon)]}
    + \norm{u_{2y}}_{L^\infty[0,T;L^2(S_\epsilon)]} 
    + \norm{\sqrt{\epsilon}(u^\epsilon_{2x} + u^\epsilon_{1y}) 
      - \gamma^\epsilon/\sqrt{\epsilon}}_{L^\infty[0,T;L^2(S_\epsilon)]} \\
    \leq C \left(
    \ltwo{u_t(0)} + \norm{u(0)}_{H^1(\Omega_0)} 
    + \norm{\gamma(0)}_{L^2(-1,1)}
    + \loneltwo{f} \right).
  \end{multline*}
  In particular, 
  $\norm{\bbC_\epsilon(\nabla
    u^\epsilon)}_{L^\infty[0,T;\Ltwo^{2 \times 2}]}$ and
  $\norm{\gammabar^\epsilon_t}_{L^2[0,T; L^2(-1,1)]}$ and
  $\norm{\gammabar^\epsilon}_{L^\infty[0,T; H^1(-1,1)]}$ are bounded where
  $\gammabar^\epsilon(t,x) = (1/\epsilon) \int_{-\epsilon/2}^{\epsilon/2}
  \gamma^\epsilon(t,x,y) \, dy$ is the average of $\gamma^\epsilon$ over
  the fault region, and the Korn and sharp Poincare inequality
  in Lemma \ref{lem:HoneEps} imply
  $$
  \norm{u^\epsilon_{1y}}_{L^\infty[0,T;L^2(S_\epsilon)]} 
  + \norm{u^\epsilon_{2x}}_{L^\infty[0,T;L^2(S_\epsilon)]} 
  \leq C/\sqrt{\epsilon}
  \quad \text{ and } \quad
  \norm{u^\epsilon}_{L^\infty[0,T;L^2(S_\epsilon)]}
  \leq C \sqrt{\epsilon}.
  $$
\end{corollary}

\subsection{Proof of Theorem \ref{thm:epsTo0}}
Fix test functions 
$$
\uhat \in \{\uhat \in H^1[0,T;H^2(\Omega_0) \times H^2(\Omega)] 
\sst \uhat(.,.,\pm 1) = 0 \text{ and } \uhat(T,.,.) = 0 \}
$$
and $\gammahat \in L^2[0,T; H^1(-1,1)]$, and note that test functions
$\uhat$ with this regularity are dense in $\{\uhat \in H^1[0,T;U_0
\times U] \sst \uhat(T) = 0\}$.  Let $\uhat^\epsilon_1 \in \Hone$ be
the function (see Figure \ref{fig:U0toU})
$$
\uhat_1^\epsilon(t,x,y) = \left\{ \begin{array}{cc}
    \uhat_1 
    \left(t,x, \frac{y-\epsilon/2}{1-\epsilon/2}\vph\right) &
    \epsilon/2 < y < 1, \\
    \left(\half + \frac{y}{\epsilon}\right) 
    \uhat_1(t,x, 0^+) +
    \left(\half - \frac{y}{\epsilon}\right)
    \uhat (t,x,0^-) &
    -\epsilon/2 \leq y \leq \epsilon/2, \\
    \uhat_1 
    \left(t,x, \frac{y+\epsilon/2}{1-\epsilon/2}\vph\right) &
    -1 < y < -\epsilon/2 ,
  \end{array} \right.
$$
and set the test functions in equations
\eqnref{:uEpsT}--\eqnref{:gEpsT} to be $\uhat^\epsilon =
(\uhat_1^\epsilon, \uhat_2) \in U_\epsilon$ and $\gammahat(t,x,y) =
\gammahat(t,x)$ and integrate the equation for $u^\epsilon$ by parts
in time to get
\begin{eqnarray}
\lefteqn{  
\int_0^T\!\int_\Omega -(\rho u^\epsilon_t, \uhat^\epsilon_t) 
+ \left( \bbC_\epsilon(D(u^\epsilon)), D(\uhat^\epsilon) \right)
+ \int_0^T\!\int_{-1}^1 \mu \left( [u^\epsilon_1] 
  - \gammabar^\epsilon,  [\uhat_1] \right) \label{eqn:uEpsA} } \\
&&
+ \int_0^T\!\int_{S_\epsilon} \mu \left(
  \epsilon (u^\epsilon_{2x} + u^\epsilon_{1y}) 
  - \gamma^\epsilon, \uhat_{2x}  \vph\right)
  + \mu( u^\epsilon_{2x} , [\uhat_1])
= \int_\Omega (\rho u^\epsilon_t(0), \uhat^\epsilon(0))
+\int_0^T\!\int_\Omega (\rho f,\uhat^\epsilon) ,
\nonumber
\end{eqnarray}
and
\begin{equation}
\int_0^T\!\int_{-1}^1 
(1/\beta) (\gammabar^\epsilon_t, \gammahat) 
+ \ell\left(\gammabar^\epsilon,\gammahat \right)
  - \mu \left([u^\epsilon_1]-\gammabar^\epsilon, \gammahat \right) 
- \int_0^T\!\int_{S_\epsilon} 
\mu (u^\epsilon_{2x}, \gammahat) = 0 ,
\label{eqn:gEpsA}
\end{equation}
where $\gammabar^\epsilon(t,x) = (1/\epsilon)
\int_{-\epsilon/2}^{\epsilon/2} \gamma^\epsilon(t,x,y)\,dy$ is the
average shear in the fault region. The last terms on the left of these
two equations represent the ``consistency error'' corresponding to
approximating a fault region of finite width with a sharp interface.
We verify that these terms vanish as $\epsilon \rightarrow 0$, and upon
passing to a sub--sequence the remaining terms consist of weakly
converging terms paired with a strongly converging test function, so
the limits of these pairings are the pairings of their
limits from which the theorem follows.

Using the bounds in Corollary \ref{cor:bounds} and Lemma \ref{lem:HoneEps}
we may pass to a subsequence for which
\begin{align*}
u^\epsilon_t &\weakst u_t, && \text{ in } \quad L^\infty[0,T;\Ltwo^2] \\
\bbC_\epsilon(D u^\epsilon) & \weakst \bbC(D u) 
&& \text{ in } \quad L^\infty[0,T;\Ltwo^{2 \times 2}] \\
[u^\epsilon_1] &\weakst [u_1], && \text{ in } \quad L^\infty[0,T;L^2(-1,1)] \\
\gammabar^\epsilon &\weakst \gamma, && \text{ in } \quad L^\infty[0,T;L^2(-1,1)] \\
\gammabar^\epsilon_t &\weak \gamma_t, && \text{ in } \quad L^2[0,T;L^2(-1,1)] \\
\gammabar^\epsilon_x &\weakst \gamma_x, && \text{ in } 
\quad L^\infty[0,T;L^2(-1,1)] 
\end{align*}
The first two terms in equations \eqnref{:uEpsA} and \eqnref{:gEpsA}
are paired with the test functions 
\begin{align*}
\uhat^\epsilon_t \rightarrow \uhat_t,& && \text{ in } \quad L^1[0,T;\Ltwo^2] \\
\begin{bmatrix} 
  \uhat^\epsilon_{1x} &  \uhat^\epsilon_{1y} \chi_{\Omega_\epsilon} \\
  \uhat_{2x} \chi_{\Omega_\epsilon} &  \uhat^\epsilon_{2y}
\end{bmatrix}
& \rightarrow \nabla \uhat
&& \text{ in } \quad L^\infty[0,T;L^2(\Omega_0)^{2 \times 2}],
\end{align*}
and the terms involving $[u^\epsilon_1]-\gammabar^\epsilon$ are
paired with test functions independent of $\epsilon$, from which it
follows that the first three terms on the left hand sides of equations
\eqnref{:uEpsA} and \eqnref{:gEpsA} converge as claimed.

The Cauchy Schwarz inequality and the smoothness of the
test function $\uhat_2(t) \in H^2(\Omega) \embed \Wonefour$ are used
to estimate the first consistency error term in equation
\eqnref{:uEpsA},
\begin{eqnarray*}
\int_0^T\!\int_{S_\epsilon} \mu \left(
  \epsilon (u^\epsilon_{2x} + u^\epsilon_{1y}) 
  - \gammabar^\epsilon, \uhat_{2x}  \vph\right)
&\leq& C \norm{\sqrt{\epsilon} (u^\epsilon_{2x} + u^\epsilon_{1y}) 
  - \gamma^\epsilon/\sqrt{\epsilon}}_{L^\infty[0,T;L^2(S_\epsilon)]}
\norm{\uhat_{2x}}_{L^1[0,T;L^2(S_\epsilon)]} \\
&\leq& C \norm{\sqrt{\epsilon} (u^\epsilon_{2x} + u^\epsilon_{1y}) 
  - \gamma^\epsilon/\sqrt{\epsilon}}_{L^\infty[0,T;L^2(S_\epsilon)]}
\norm{\uhat_{2x}}_{L^1[0,T;L^4(S_\epsilon)]} \epsilon^{1/4} \\
&\rightarrow& 0.
\end{eqnarray*}
The final terms on the left hand side of equations \eqnref{:uEpsA} and
\eqnref{:gEpsA} involve $u^\epsilon_{2x}$ paired with test functions
which are independent of $y$. It then suffices to show that
$\ubar^\epsilon_{2x} \equiv \int_{-\epsilon/2}^{\epsilon/2}
u^\epsilon_{2x}(.,.,y) \, dy$ converges weakly star to zero in
$L^\infty[0,T; L^2(-1,1)]$.  To do this the Cauchy Schwarz inequality
and Corollary \ref{cor:bounds} are used to first show that it is bounded,
$$
%\norm{\ubar^\epsilon_{2x}(t)}_{L^2(-1,1)}^2 = 
\int_{-1}^1 \left( 
  \int_{-\epsilon/2}^{\epsilon/2} u^\epsilon_{2x}(t,x,y) \, dy 
\right)^2 \, dx
\leq \int_{-1}^1 \epsilon
  \int_{-\epsilon/2}^{\epsilon/2} u^\epsilon_{2x}(t,x,y)^2 \, dy \, dx
= \epsilon \norm{u^\epsilon_{2x}(t)}_{L^2(S_\epsilon)}^2
\leq C.
$$
To establish weak star convergence to zero it then suffices to
test against smooth functions $\phihat$ with compact support in $(0,T)
\times \Omega$ since they are dense in $L^1[0,T;L^2(-1,1)]$,
$$
\int_0^T \! \int_{-1}^1 \! \int_{-\epsilon/2}^{\epsilon/2} 
(\uhat^\epsilon_{2x}, \phihat)
= \int_0^T \! \int_{-1}^1 \! \int_{-\epsilon/2}^{\epsilon/2} 
-(\uhat^\epsilon_2, \phi_x)
\leq \norm{\uhat^\epsilon_2}_{L^\infty[0,T;L^2(S_\epsilon)]} 
 \norm{\phihat}_{L^1[0,T;L^2(S_\epsilon)]} 
\leq C \sqrt{\epsilon}.
$$
It follows that the limit $(u,\gamma)$ is a solution of the sharp interface 
problem, and the theorem follows provided it takes the specified initial values.
However, this is direct since $(u^\epsilon, \gammabar^\epsilon)$ converges
weakly in $H^1[0,T; \Ltwo^2 \times L^2(-1,1)]$ from which it follows that
the initial values of the limit $(u,\gamma)$ are the limit of the initial
values.

\section{Numerical Examples} \label{sec:NumericalExamples}
This section first presents a numerical example to exhibit the contrast
between direct numerical simulation of the stationary form of equations
\eqnref{:momentum}--\eqnref{:gamma} and the limit problem for the
considered in Section 2. In the second section numerical approximation
of a singular solution corresponding to a dislocation is presented.

In the numerical examples below the parameters are set to
$$
\mu=1, \quad \lambda = 2, \quad a=1/2, \quad \ell = 1/4, \quad 
\etahat=2, \quad \nu = 0, \quad \epsilon = 1/10.
$$
and for the limit problem uniform rectangular elements of size $h =
1/n$ with $n \in \bbN$ are utilized. When $\epsilon > 0$ the fault
region is meshed with rectangular elements of size $1/n \times
\epsilon/n$; the mesh with $n = 4$ is illustrated in Figure
\ref{fig:elastic3m844}. Galerkin approximations of the solution to the
elasticity problems are computed using the piecewise quadratic finite
element spaces on these meshes.

\subsection{Classical Solution}
A piecewise smooth solution of the limit problem with constant Lame
parameters is constructed by setting
\begin{equation} \label{eqn:u0}
u(x,y) = \left\{ \begin{array}{ll}
\frac{1}{2} \left( \begin{array}{c} e^{a(x-y)} \\ 
e^{a(x-\kappa y)} \end{array} \right) 
+ \left( \begin{array}{r} \phi_y(x,y) \\ 
-\phi_x(x,y) \end{array} \right) \qquad
& y > 0, \\ [3ex]
\frac{1}{2} \left( \begin{array}{c} -e^{a(x+y)} \\ 
e^{a(x+\kappa y)} \end{array} \right) 
+ \left( \begin{array}{r} \phi_y(x,y) \\ 
-\phi_x(x,y) \end{array} \right)
& y < 0,
\end{array} \right.
\end{equation}
where $\kappa = \lambda/(2\mu + \lambda)$ and 
$\phi(x,y) = e^{-\ell y}\cos(\ell x)$. Then
\begin{equation} \label{eqn:g0}
[u_1(x)] = e^{ax},
\qquad
\gamma(x) \equiv [u_1] - (1/\mu) \bbC(\nabla u)_{12}
= e^{ax} - 2 \ell^2 \cos(\ell x),
\end{equation}
and right hand sides for the stationary problem are manufactured so
that the equations are satisfied,
$$
f = -\dv(\bbC(\nabla u))
\qquad \text{ and } \qquad
f_0 = \etahat \gamma - \bbC(\nabla u)_{12} - \nu \gamma_{xx}.
$$

To exhibit the differences between direct numerical approximation of
\eqnref{:momentum}--\eqnref{:gamma} and numerical approximation of the
limiting problem we first tabulate the errors, $u^0 - u^0_h$, of the
numerical approximation of the solution \eqnref{:u0}.  Numerical
approximations $u^\epsilon_h$ of the stationary equations
\eqnref{:momentum}--\eqnref{:gamma} are then computed using the same
boundary data and body force $f$. While the exact solution,
$u^\epsilon$, of the problem with this data is not known, we tabulate
(norms of) the differences $u^0 - u^\epsilon_h$ for $\epsilon$
fixed. As $h \rightarrow 0$ this difference converges to the
``modeling'' error $u^0 - u^\epsilon$ associated with approximating
the fault region by a surface. An estimate of the mesh
size required to resolve the deformation in the fault region is
obtained by observing when difference $u^0 - u^\epsilon_h$ stabilizes.
Note that in general $\lim_{\epsilon \rightarrow} \ltwo{u^0 -
  u^\epsilon} \rightarrow 0$ but $\lim_{\epsilon \rightarrow}
\norm{u^0 - u^\epsilon}_{H^1(\Omega_0)} \not\rightarrow 0$.

\subsubsection{Uncoupled Problem}
\begin{table}[h]
\begin{center}
  \begin{tabular}{ |c | c | c | r|}
    \hline
    $h$ & $\ltwo{u^0-u^0_h}$   & $\norm{u^0-u^0_h}_{H^1(\Omega_0)}$ & 
    \# unknowns \\ \hline
    $1/8$ & 1.684735e-05  &  4.402589e-04 & 629 \\ 
    $1/16$ & 2.098734e-06 & 1.091350e-04 & 2277 \\ 
    $1/32$ & 2.617986e-07 & 2.718018e-05 & 8645 \\ 
    $1/64$ & 3.268957e-08 & 6.783206e-06 & 33669 \\ 
    $1/128$ & 4.084190e-09  & 1.694404e-06 & 132869 \\ 
    $1/256$ & 5.154866e-10 & 4.234319e-07&  527877 \\  \hline
    Norms & 1.435134 & 1.662724 & \\
    \hline
  \end{tabular} 
  \caption{Errors for the uncoupled limit problem ($\epsilon = 0$).}
  \label{uncoupled_limit_table}
\end{center}
\end{table}
\begin{table}[h]
\begin{center}
  \begin{tabular}{ | c | c | c | r |}
    \hline
    $h$ & $\|u^0-u^\epsilon_h\|_\Ltwo$  
    & $\|u^0-u^\epsilon_h\|_{H^1(\Omega_0)}$ & \# unknowns \\ \hline
     $1/8$   & 2.148338e-01 &  4.553696e+00 &1666 \\ 
     $1/16$  & 3.029500e-01 & 5.114175e+00 & 6402 \\ 
     $1/32$  & 3.183221e-01  & 5.119497e+00 & 25090 \\
     $1/64$  & 3.358588e-01  & 5.177173e+00 & 99330 \\ 
     $1/128$ & 3.459277e-01  & 5.223399e+00 & 395266 \\ 
     $1/256$ & 3.433117e-01  & 5.218226e+00 & 1576962	\\
    \hline
    \end{tabular}
    \caption{Differences between the uncoupled
        $\epsilon$-problem and limit problem with $\epsilon=0.1$.}
    \label{uncoupled_eps0.1}
    \end{center}
\end{table}
\begin{figure}[H]
\begin{center}
\includegraphics[width=4.0in]{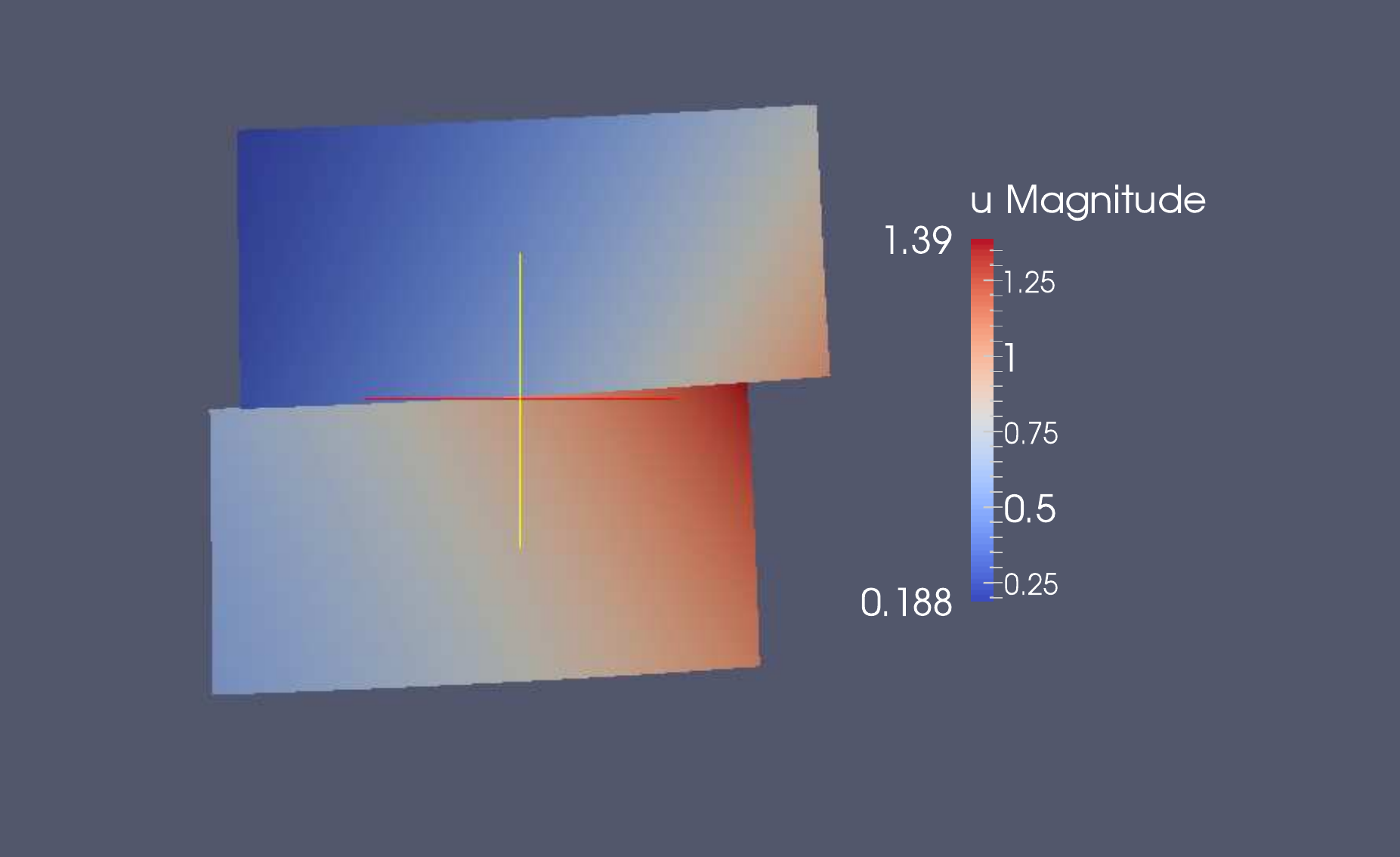}
%\vspace*{-1.75in}
\caption{Solution of the uncoupled limit problem mesh size $h=1/256$.}
\label{fig:uncoupled_limit_problem}
\end{center}
\end{figure}
\begin{figure}[H]
\begin{center}
\includegraphics[width=4.0in]{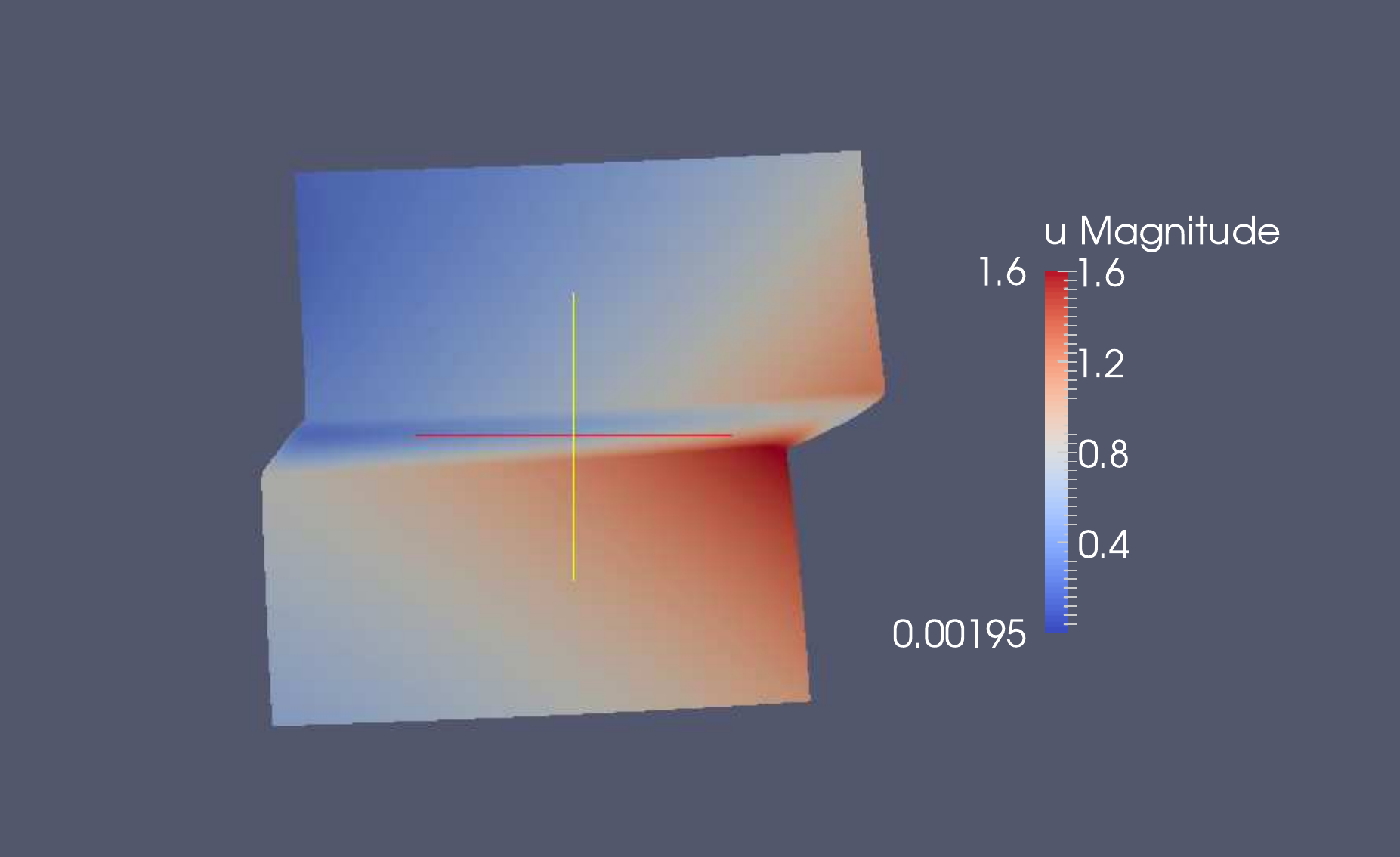}
\caption{Solution of the uncoupled problem with 
    $\epsilon=0.1$, $h=1/256$.}
\label{fig:uncoupled_256_eps_zero_pt1}
\end{center}
\end{figure}

Table \ref{uncoupled_limit_table} exhibits the errors in the numerical
approximation of the solution \eqnref{:u0} of the limit problem
considered in Section \ref{sec:Gamma} with $\gamma$ the function
specified in \eqnref{:g0}. The optimal third order rate in $\Ltwo$ and
second order rate for the derivatives is clear.  Norms of the
differences $u^0 - u^\epsilon_h$ are presented in Table
\ref{uncoupled_eps0.1}. For this example it is clear that very
accurate solutions of the limit problem can be computed on very modest
meshes while resolution of the deformation in the fault region
requires significantly finer meshes. The norms computed on the finest
meshes give an estimate of the modeling error $\ltwo{u^0 - u^\epsilon}
\simeq 0.34$. Representative solutions for each of the problems are
illustrated in Figures \ref{fig:uncoupled_limit_problem} and
\ref{fig:uncoupled_256_eps_zero_pt1}.
 
\subsubsection{Coupled Problems}
\begin{table}[H]
  \begin{center}
    \begin{tabular}{ | c | c | c | c | r |} \hline
      $h$ & $\ltwo{u^0-u^0_h}$   & $\norm{u^0-u^0_h}_{H^1(\Omega_0)}$ & 
      $\norm{\gamma-\gamma_h}_{L^2(-1,1)}$ & \# unknowns \\ \hline
     $1/8$ & 1.683317e-05  &  4.403035e-04 & 1.597969e-05& 646\\
     $1/16$ & 2.098047e-06 & 1.091376e-04 & 2.077018e-06& 2310\\
     $1/32$ & 2.617686e-07 & 2.718031e-05 & 2.646097e-07& 8710\\
     $1/64$ & 3.268836e-08 & 6.783211e-06 & 3.338419e-08& 33798\\ 
     $1/128$& 4.084169e-09  & 1.694404e-06 & 4.191933e-09& 133126\\
     $1/256$& 5.161811e-10  & 4.234319e-07& 5.258066e-10 & 528390\\ \hline
     Norms & 1.435134 & 1.662724 & 1.365834 & \\
     \hline
   \end{tabular} 
    \caption{Errors for the coupled limit problem ($\epsilon = 0$).}
    \label{limit_table}
  \end{center}
\end{table}    
\begin{table}[H]
  \begin{center}
    \begin{tabular}{ | c | c | c | c | r |}
      \hline
      $h$  & $\ltwo{u^0-u^\epsilon_h}$ & 
      $\norm{u^0-u^\epsilon_h}_{H^1(\Omega_0)}$  & 
      $\norm{\gamma - \gamma^\epsilon_h}_{L^2(S_\epsilon)}$ & 
      \# unknowns \\ \hline
      $1/8$ & 1.575165e-01     &  4.626232e+00 & 4.662827e-02 & 1955\\
      $1/16$ & 1.582971e-01   &  4.628384e+00 & 5.548420e-02 & 7491 \\
      $1/32$ & 1.584715e-01   & 4.655150e+00 &6.875543e-02 & 29315\\ 
      $1/64$ & 1.594472e-01   & 4.670685e+00 & 7.885558e-02 & 115971\\
      $1/128$& 1.596179e-01   & 4.675027e+00 & 8.156652e-02 & 461315 \\
      $1/256$& 1.593571e-01   & 4.676821e+00 & 8.152622e-02 & 1840131\\
      \hline
    \end{tabular}
    \caption{Differences between coupled $\epsilon$-problem and
        limit problem with $\epsilon=0.1$.}
    \label{eps0.1table}
  \end{center}
\end{table}
\begin{figure}[H]
  \begin{center}
    \includegraphics[width=4.0in]{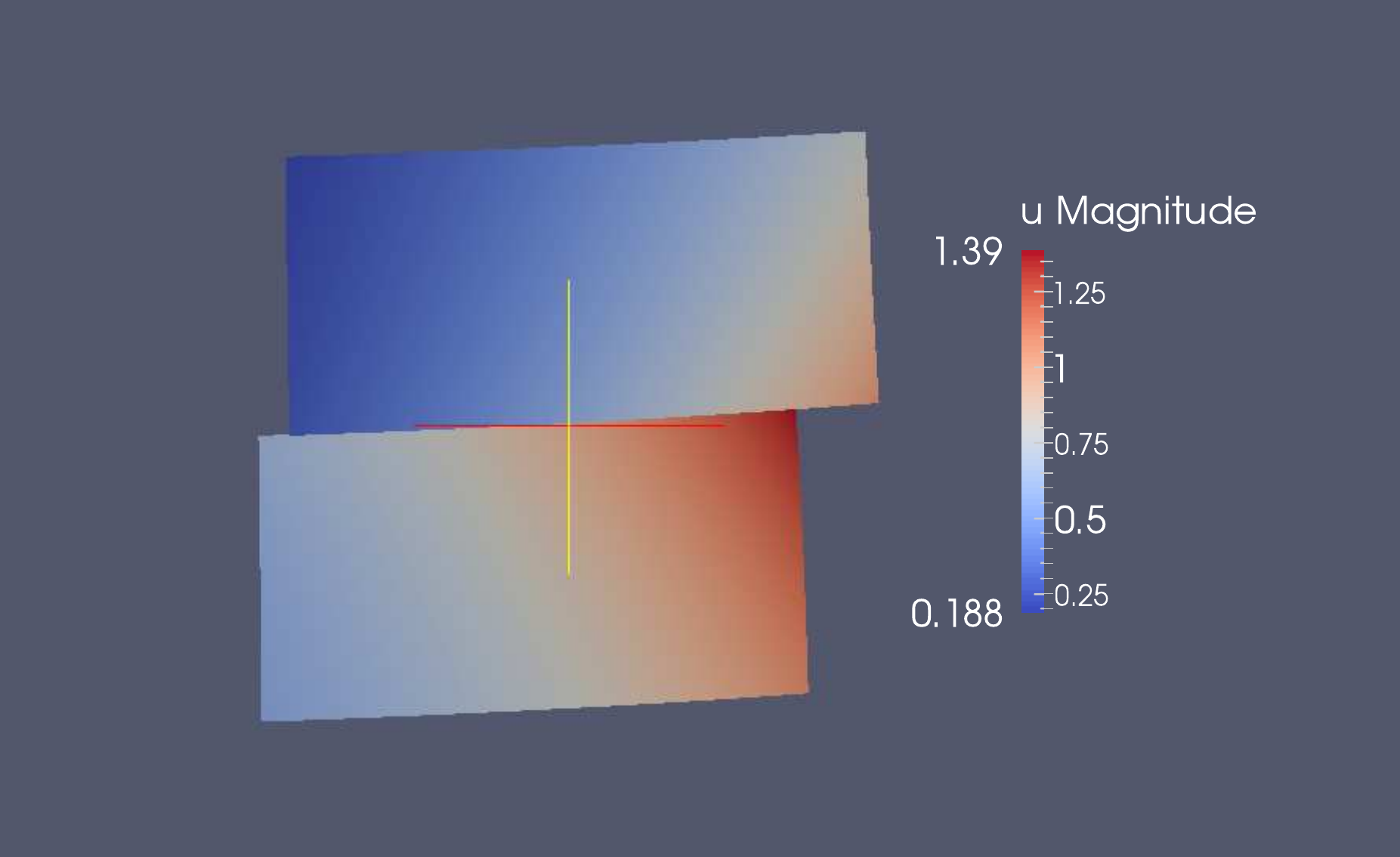}
    \caption{Displacement of coupled limit problem with $h=1/256$.}
    \label{fig:elastic0_256}
  \end{center}
\end{figure}
\begin{figure}[H]
\begin{center}
\includegraphics[width=4.0in]{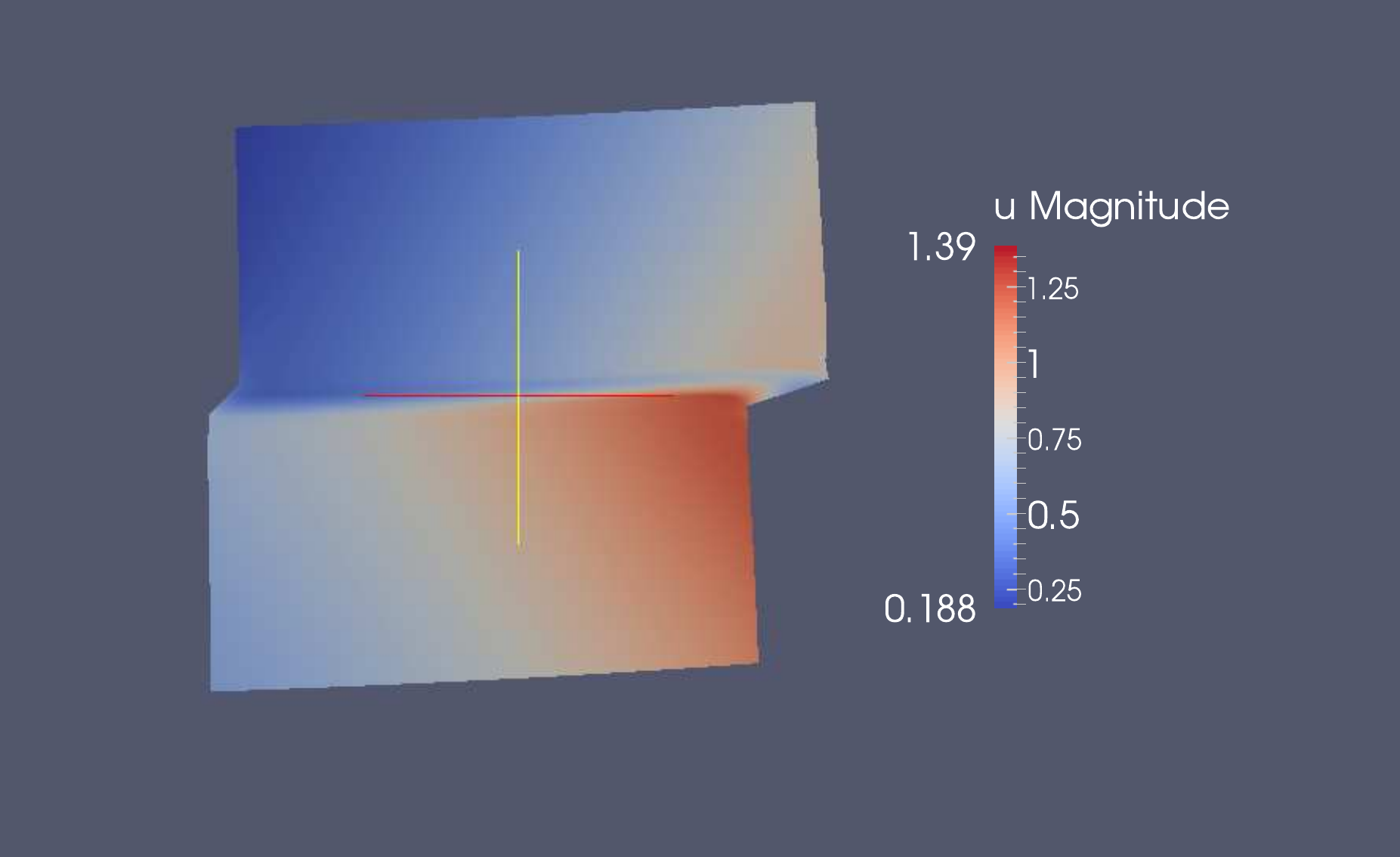}
\caption{Displacement for coupled problem with $\epsilon=0.1$ and $h=1/256$.}
\label{fig:elastic_256_0pt1}
\end{center}
\end{figure}

Table \ref{limit_table} exhibits the errors for the coupled problem
when numerical approximations of both $u$ and $\gamma$ are computed
using the limit energy given in equation \eqnref{:coupled}. Again the
optimal third order rate in the $L^2$ norms for both $u$ and $\gamma$
and second order rate for the derivatives of $u$ is obtained.  Norms
of the differences $u^0 - u^\epsilon_h$ and $\gamma -
\gamma^\epsilon_h$ are presented in Table \ref{eps0.1table}. As for
the uncoupled case, very accurate solutions of the limit problem can
be computed on modest meshes while resolution of the deformation in
the fault region requires finer meshes.  The modeling errors for this
problem are $\ltwo{u^0 - u^\epsilon} \simeq 0.16$ and
$\norm{\gamma-\gamma^\epsilon} \simeq 0.082$. Representative
deformations are illustrated in Figures \ref{fig:elastic0_256} and
\ref{fig:elastic_256_0pt1}.

\subsection{Dislocation}
\begin{figure}
\begin{center}
\includegraphics[width=2.25in,clip=true, bb=55 435 370 740]{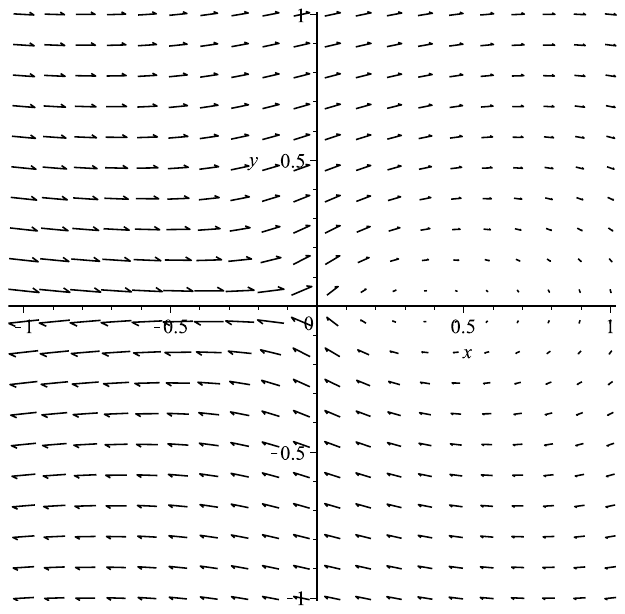}
\qquad
\includegraphics[width=2.25in,clip=true, bb=95 350 500 740]{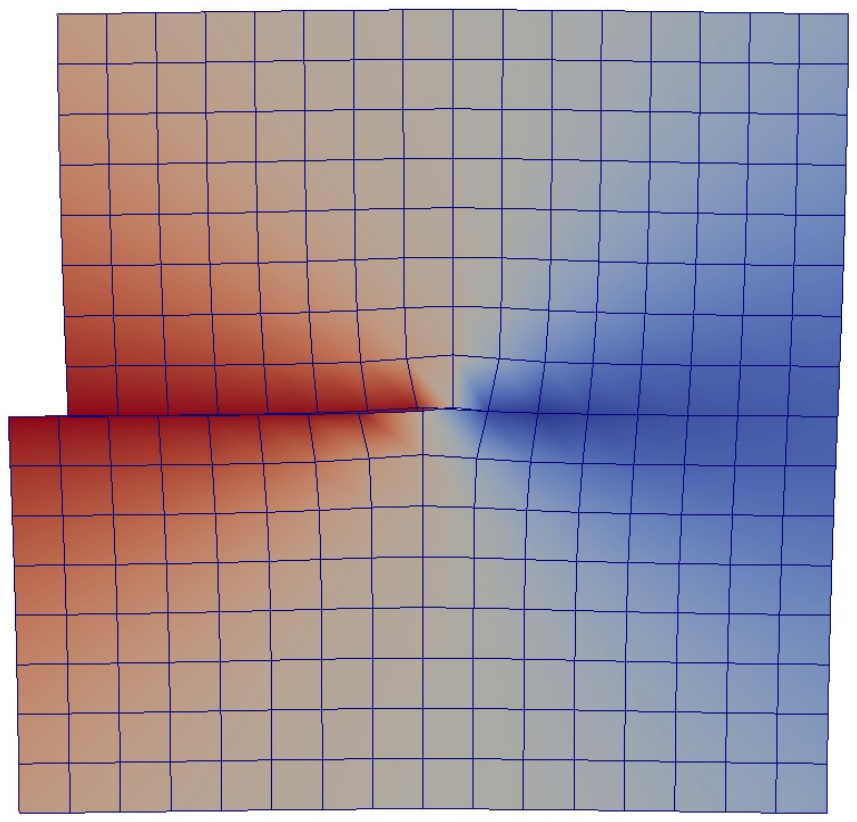}
\caption{Displacement field with dislocation at the origin and numerical approximation.}
\label{fig:dislocationU}
\end{center}
\end{figure}

An explicit solution for the linear elasticity problem
with an edge dislocation along the $z$--axis is \cite{HiLo82}
$$
u(x,y) = \frac{1}{2\pi} \begin{pmatrix}
\arctan(y/x) + \frac{xy}{2 (1-\nuhat)(x^2+y^2)} \\
\frac{-1}{4(1-\nuhat)} \left( (1-2\nuhat) \ln(x^2+y^2) 
+ \frac{x^2-y^2}{x^2+y^2} \right)
\end{pmatrix},
\text{ with Poisson ratio } \nuhat = \frac{\lambda}{2(\mu+\lambda)}.
$$
This solution, illustrated in Figure \ref{fig:dislocationU},
represents the displacement that results when a dislocation, currently
at the origin, has propagated along the negative $x$--axis so that
$[u_1(x)] = 1$ for $x < 0$ and $[u_1(x)] = 0$ for $x > 0$. The stress
has a singularity of order $O(1/r)$ at the origin and is otherwise
continuous, and the displacement is square integrable but its derivatives
are not. While the results of the prior sections are not
applicable to singular solutions, almost singular solutions arise in
engineering practice so it is important for the numerical schemes 
to be robust in this context. 

\begin{table}[H]
  \begin{center}
    \begin{tabular}{|c|c|c|}
      \hline
      $h$  & $\ltwo{u-u_h}$ & $\norm{\gamma_h}_{L^2(-1,1)}$ \\ \hline
      $1/8$ & 1.226581e-02 &  7.005017 \\
      $1/16$ & 6.141691e-03 & 9.835010 \\
      $1/32$ & 3.073164e-03 & 13.85832\\
      $1/64$ & 1.537180e-03 & 19.56296\\
      $1/128$& 7.687419e-04 & 27.64101 \\
      $1/256$& 3.844092e-04 & 39.07248 \\
      \hline
    \end{tabular}
    \caption{Errors in displacement and norm of shear for 
      dislocation example.}
    \label{tbl:dislocation}
  \end{center}
\end{table}

To illustrate the robustness properties of codes using the limit
energy a singular solution of the stationary limit problem is
manufactured by setting $\gamma = [u_1] - (1/\mu) T_{12}$, so that the
jump condition is satisfied, and non--homogeneous right hand side for
the equation for $\gamma$,
$$
f_0(x) = \etahat \gamma(x) - T_{12}(x,0) 
= \etahat \, [u_1(x)]
+ \frac{(\lambda+\mu)(\etahat+\mu)}{\pi (2\mu+\lambda) x}.
$$
(Since $\gamma_{xx}$ does not exist we set the coefficient of this
term to be zero.)  Inner products of this (non--integrable) function
with basis functions were approximated using Gaussian quadrature.
Table \ref{tbl:dislocation} shows that the error $\ltwo{u-u_h}$
converges linearly with $h$ and $\ltwo{\gamma_h} \simeq O(1/\sqrt{h})$
diverges since the limit $\gamma(x) \simeq O(1/x)$ is not
integrable. Figure \ref{fig:dislocationU} illustrates the deformation
computed with quadratic elements on a $16 \times 16$ grid.

\appendix
\section{Derivation from a Plasticity Model}
Displacements and gradients are assumed to be small in the region
$\Omega_\epsilon = \Omega \setminus \Sbar_\epsilon$ outside
the fault so that the motion is governed by the equations
of linear elasticity,
$$
\rho u_{tt} - \dv( \bbC(\nabla u) ) = \rho f,
\qquad \text{ in } \Omega_\epsilon.
$$
Small displacement 
plasticity theory models the motion in the fault region $S_\epsilon$. In this theory
the elastic deformation tensor, $U$, deviates from $\nabla u$ due 
to slips and motion of defects. The balance of linear momentum becomes
$$
\rho u_{tt} - \dv( \bbC(U) ) = \rho f,
\qquad \text{ in } S_\epsilon,
$$
and evolution of $U$ is governed by an equation of the form
\begin{equation} \label{eqn:Ueqn}
U_t - \nabla u_t + \Curl(U) \times v_d = 0,
\end{equation}
where $v_d$ is a constitutively specified defect velocity. In this
equation the $\Curl(.)$ and cross product of a matrix act row--wise;
$$
\Curl(U)_{mn} = \epsilon_{ijn} U_{mj,i},
\qquad \text{ and } \qquad
(A \times v)_{mn} = \epsilon_{ijn} A_{mi} v_j.
$$
The defect velocity $v_d$ is chosen to model the (typically large)
dissipation due to defect motion, and local energy changes due to
distortion in the material during passage of a defect. In the
following lemma the axial vector of the skew part of a matrix $A$ is
denoted by $X(A)$; that is
$$
X(A)_i = \epsilon_{ijk} A_{jk}.
$$

\begin{lemma}
  Let $\eta:\Re \rightarrow \Re$ and $\beta,\gamma,T_{12}:(0,T) \times
  S_\epsilon \rightarrow \Re$ be smooth, $\nu \in \Re$, and suppose
  that
  $$
  \gamma_t 
  + \beta \gamma_x^2 \left(\eta'(\gamma) 
    - \nu \gamma_{xx} - T_{12} \vph\right) = 0,
  \qquad \text{ on } (0,T) \times S_\epsilon.
  $$
  Let
  $$
  u:(0,T) \times S_\epsilon \rightarrow \Re^2 \embed \Re^3,
  \qquad \text{ and } \qquad
  T:(0,T) \times S_\epsilon \rightarrow 
  \Re^{2 \times 2}_{sym} \embed \Re^{3 \times 3}_{sym},
  $$
  be smooth with $T_{12}$ as above, and let
  $$
  U = \nabla u -
  \begin{bmatrix} 
    0 & \gamma & 0 \\ 0 & 0 & 0 \\ 0 & 0 & 0
  \end{bmatrix}
  \qquad \text{ and } \qquad
  v_d = 
  \left(I - \frac{\omega}{|\omega|} \otimes \frac{\omega}{|\omega|} \right) 
  X\left(S \, \Curl(U) \vph\right),
  $$
  where
  $$
  S = \left(T -
    \begin{bmatrix} 0 & 2 \eta'(\gamma) & 0 
      \\ 0 & 0 & 0 \\ 0 & 0 & 0 \end{bmatrix} 
    + 2 \nu \Curl(\Curl(U)) \right)_{sym}
  \quad \text{ and } \quad
  \omega = X(\Curl(U)).
  $$
  Then the triple $(U,u,v_d)$ satisfies equation \eqnref{:Ueqn}.
\end{lemma}

Under the ansatz of the lemma the matrices
$\Curl(U)$ and $\Curl(\Curl(U))$ become
$$
\Curl(U) = \begin{bmatrix} 
    0 & 0 & -\gamma_x \\ 0 & 0 & 0 \\ 0 & 0 & 0
  \end{bmatrix},
\qquad
\Curl(\Curl(U)) = \begin{bmatrix} 
    0 & -\gamma_{xy} & \gamma_{xx} \\ 0 & 0 & 0 \\ 0 & 0 & 0
  \end{bmatrix}
$$
and vectors $X(\Curl(U))$ and $v_d$ are
$$
X(\Curl(U)) = \begin{pmatrix} 0 \\ -\gamma_x \\ 0 \end{pmatrix},
\qquad
v_d = \begin{pmatrix} \beta \gamma_x \left( \eta'(\gamma) 
    - \nu \gamma_{xx} - T_{12} \vph\right) \\ 0 \\ 0 \end{pmatrix}.
$$

{\bf Acknowledgement:} The authors appreciate the insight gained from
multiple discussions with Professor Amit Acharya of the Department of
Civil and Environmental Engineering at Carnegie Mellon University.

\bibliography{references}

\end{document}